\documentclass[11pt,a4paper]{article}
\usepackage{amsfonts,amsgen,amstext,amsbsy,amsopn,amsfonts,amssymb,amscd}
\usepackage[leqno]{amsmath}
\usepackage[amsmath,amsthm,thmmarks]{ntheorem}
\usepackage{float}
\usepackage{color}
\usepackage{tikz}
\usepackage{multirow}
\usepackage{mathrsfs}
\usepackage{graphicx}
\usepackage{subfigure}
\newtheorem{thm}{Theorem}[section]

\newtheorem{prop}[thm]{Proposition}

\newtheorem{lem}[thm]{Lemma}

\newtheorem{false statement}{False statement}

\theoremstyle{definition}

\makeatletter \@addtoreset{equation}{section}

\def\hg{\mathcal{G}}

\begin{document}

\title{\bf\Large Large Maximal Intersecting Families with Covering Number at Least Three}
\date{}
\author{ Huajun Zhang\\[10pt]
Shaoxing University,
Shaoxing 312000 P. R. China\\[6pt]
 huajunzhang@usx.edu.cn}
\maketitle
\begin{abstract}
For all $n>2k$ and $k\ge4$, we determine the second largest size
of a maximal intersecting family $\mathcal F\subseteq\binom{[n]}k$
whose covering number $\tau(\mathcal F)$ is at least three.
Apart from the largest family $\mathcal G_1(n,k)$, the sharp bound is
$\max\{q(n,k),|\mathcal G_2(n,k)|\}$, where $q(n,k)$ is attained by
families with three exceptional members. For $k\ge5$, the transition
between the two bounds occurs at $n=3k-2$; for $k=4$, they coincide.
We also classify the equality cases under a critical shifting
hypothesis, including additional configurations at $n=2k+1$.
The proof combines exact sections relative to a minimum three-cover
with local replacements preserving the covering number.
A quantitative inverse-shift bound, obtained from a disjointness
graph, controls maximal completions and yields the upper bound
without a shifting assumption on the original family.
\end{abstract}

\section{Introduction}

Let $[n]=\{1,2,\ldots,n\}$.
We interpret $\binom ab$ as zero if $a<0$, $b<0$, or $b>a$.
Two families are isomorphic, denoted by $\cong$, if a permutation of
the ground set maps one onto the other. For $1\le i<j\le n$, write
$[i,j]=\{i,i+1,\ldots,j\}$. Let $2^{[n]}$ denote the power set of $[n]$,
and let $\binom{[n]}{k}$ denote the family of all $k$-element subsets of
$[n]$.
A family $\mathcal F\subseteq \binom{[n]}{k}$ is called
\emph{intersecting} if
$
F\cap F'\neq\emptyset
\quad\text{for all } F,F'\in\mathcal F.
$
It is called \emph{trivial} if
$
\bigcap_{F\in\mathcal F}F\neq\emptyset,
$
and \emph{nontrivial} otherwise. Similarly, two families $\mathcal A$ and
$\mathcal B$ are called \emph{cross-intersecting} if
$
A\cap B\neq\emptyset$
\text{\ for all } $A\in\mathcal A$ \text{ and } $B\in\mathcal B$.

A family $\mathcal F \subseteq \binom{[n]}{k}$ is called \emph{maximum} (with respect to a given class) if it attains the largest possible cardinality among all families in that class.
On the other hand, $\mathcal F$ is called \emph{maximal} if it is not properly contained in any other family of the same class.
Clearly, every maximum family is maximal, but the converse does not hold in general.
A set $X\subseteq [n]$ is a \emph{cover}, or \emph{transversal}, of
$\mathcal F$ if
$
X\cap F\neq\emptyset$ for all $F\in\mathcal F.$
The \emph{covering number} of $\mathcal F$, denoted by $\tau(\mathcal F)$,
is the minimum size of a cover of $\mathcal F$. Thus, a nonempty
intersecting family $\mathcal F$ is nontrivial if and only if
$\tau(\mathcal F)\ge 2$.

One of the central results in extremal set theory is the Erd\H{o}s--Ko--Rado theorem.

\begin{thm}[Erd\H{o}s--Ko--Rado \cite{EKR}]
Let $\mathcal F\subseteq\binom{[n]}{k}$ be intersecting with $n\ge2k$. Then
\[
|\mathcal F|\le\binom{n-1}{k-1}.
\]
Moreover, if $n>2k$, equality holds if and only if
$
\mathcal F=\{F\in\binom{[n]}{k}:x\in F\}
$
for some $x\in[n]$.
\end{thm}

Hence, for $n\ge2k+1$, the largest intersecting families are trivial.
For nontrivial intersecting families, Hilton and Milner proved the following result.

\begin{thm}[Hilton--Milner \cite{HM}]
Let $n\ge2k+1$.
If $\mathcal F\subseteq\binom{[n]}{k}$ is nontrivial and intersecting, then
\[
|\mathcal F|
\le
\binom{n-1}{k-1}
-
\binom{n-k-1}{k-1}
+1.
\]
\end{thm}

For $1\le t\le n-k$, define
\[
\mathcal H(n,k;t)
=
\left\{
A\in\binom{[n]}{k}:
1\in A \text{ and }
A\cap B\neq\emptyset \text{ for all } B\in\mathcal M
\right\}
\cup
\mathcal M,
\]
where
\[
\mathcal M=\{[2,k]\cup\{j\}:k+1\le j\le k+t\}.
\]

Then $\mathcal H(n,k;1)$ attains the Hilton--Milner bound and is therefore the largest nontrivial intersecting family.
In particular, it is the largest intersecting family with covering number at least $2$.
Moreover, $\tau(\mathcal H(n,k;t))=2$ whenever $k+t\le n$.
For $k\ge4$ and $n>2k$, Han and Kohayakawa \cite{hk2017} showed that
$\mathcal H(n,k;2)$ attains the maximum size among nontrivial
intersecting families not contained in any copy of the Hilton--Milner
family, and characterized all equality cases. Kostochka and Mubayi
\cite{Kostochka-Mubayi2017} investigated the families
$\mathcal H(n,k;t)$ for $t\ge3$. More recently,
Kupavskii \cite{Kupavskii2025} investigated the structure of large non-trivial
intersecting families and obtained a precise description of families
close to the Hilton--Milner bound. These results further highlight the
importance of structural stability beyond determining only the maximum
size. A natural measure of how far the intersecting family is from the
trivial family is the covering number.

Define
\[
f(n,k,s)
=
\max\bigl\{
|\mathcal F|:
\mathcal F\subseteq\binom{[n]}{k} \text{ is intersecting and }
\tau(\mathcal F)\ge s
\bigr\}.
\]
With this notation, the Erd\H{o}s--Ko--Rado and Hilton--Milner theorems state that
\[
f(n,k,1)=\binom{n-1}{k-1},
\qquad
f(n,k,2)
=
\binom{n-1}{k-1}
-
\binom{n-k-1}{k-1}
+1
\]
for $n\ge2k+1$.

\subsection{The second largest size}

We begin with the constructions associated with covering number three.
For $1\le p\le k-1$, define
\[
\mathcal B_p
=
\{[2,k]\cup\{j\}:k+1\le j\le k+p+1\}
\cup
\{[k+1,2k]\},
\]
\[
\mathcal A_p
=
\{A\in\binom{[n]}{k}:
1\in A \text{ and }
A\cap B\neq\emptyset \text{ for each } B\in\mathcal B_p
\},
\]
and
\[
\mathcal G_p(n,k)=\mathcal A_p\cup\mathcal B_p.
\]
Clearly, $\mathcal G_p(n,k)$ is maximal intersecting with covering
number three. Inclusion--exclusion gives
\[
|\mathcal G_p(n,k)|
=
\binom{n-1}{k-1}
-
2\binom{n-k-1}{k-1}
-
\sum_{j=2}^{p+1}\binom{n-k-j}{k-j}
+
\binom{n-2k}{k-1}
+
p+2.
\]
Frankl \cite{F80} first proved that for $k \geq 4$ and $n > n_0(k)$,
\begin{equation}\label{fw1}
    f(n,k,3) = |\hg_1(n,k)|,
\end{equation}
and that $\hg_1(n,k)$ is the unique optimal family up to isomorphism. Recently, Frankl and Wang \cite{FW} proved that \eqref{fw1} holds for $n > 2k$ and $k \geq 7$, while Kupavskii \cite{AK} independently established the result for $k \geq 100$. Subsequently, Frankl and Wang \cite{FW1} settled the remaining cases $k \in \{4,5,6\}$, thereby completing the proof for all $k \geq 4$ and $n > 2k$.

Our aim is to determine the next possible cardinality among maximal
intersecting families with covering number at least three, throughout
the full range $n>2k$ and $k\ge4$.

Write $g_p(n,k)=|\mathcal G_p(n,k)|$, and call the distinguished
point $1$ the \emph{hub} of these constructions. Put
\begin{equation}\label{eq:Gp-difference}
g_p(n,k)-g_{p+1}(n,k)
=\binom{n-k-p-2}{k-p-2}-1\qquad(1\le p\le k-2).
\end{equation}
In particular, $g_1>g_2\ge g_p$ for $2\le p\le k-1$, with
$g_p=g_2$ exactly when $p=2$, or when $k=4$ and $p=3$.

Two constructions with three exceptional members compete with
$\mathcal G_2(n,k)$ for the second largest size. Let
\[
B_1=[2,k+1],
\qquad
B_2=[2,k-1]\cup\{k+2,k+3\},
\]
and put
\[
B_3^{(1)}=[k,2k]\setminus\{k+3\},
\qquad
B_3^{(2)}=[k,2k-1].
\]
For $j\in\{1,2\}$, define
\begin{align*}
\mathcal Q^j(n,k)
&=\{B_1,B_2,B_3^{(j)}\}\\
&\quad\cup
\left\{
F\in\binom{[n]}k:
1\in F,\;
F\cap B\ne\emptyset\text{ for every }B\in\{B_1,B_2,B_3^{(j)}\}
\right\}.
\end{align*}
Both families are clearly maximal intersecting with covering number
three and have the same size. Inclusion--exclusion gives,
for $j\in\{1,2\}$,
\begin{align*}
q(n,k):=|\mathcal Q^j(n,k)|
={}&
\binom{n-1}{k-1}
-3\binom{n-k-1}{k-1}+\binom{n-k-3}{k-1}
+\binom{n-2k+1}{k-1}+3.
\end{align*}
Set
$M(n,k)=\max\{q(n,k),g_2(n,k)\}$.

\begin{thm}[Second largest size]\label{thm:second-size}
Let $n>2k$ and $k\ge4$.
Suppose that $\mathcal F\subseteq\binom{[n]}k$ is a maximal
intersecting family with $\tau(\mathcal F)\ge3$.
Then either
$\mathcal F\cong\mathcal G_1(n,k)$ or
\[
|\mathcal F|\le M(n,k).
\]
The bound is sharp:
$\mathcal Q^1(n,k)$ and $\mathcal Q^2(n,k)$ attain it whenever
$q(n,k)\ge g_2(n,k)$, while $\mathcal G_2(n,k)$ attains it
whenever $g_2(n,k)\ge q(n,k)$.
\end{thm}

For $k=4$, the two terms coincide:
$q(n,4)=g_2(n,4)=12n-61, \quad g_1(n,4)=13n-69$.
For $k\ge5$, their relative sizes depend on $n$:
\[
\begin{cases}
q(n,k)>g_2(n,k),
& 2k+1\le n\le3k-3,\\
q(n,k)=g_2(n,k),
& n=3k-2,\\
q(n,k)<g_2(n,k),
& n\ge3k-1.
\end{cases}
\]
In particular, $M(n,k)<g_1(n,k)$, so the theorem identifies the exact
gap below the unique largest family. It also determines where the
constructions attaining the second largest size change: the threshold
$n=3k-2$ separates the two regimes for $k\ge5$, while $k=4$ gives
equal construction sizes throughout. These comparisons follow directly
from the displayed formulas for $g_1$, $g_2$, and $q$.
The numerical bound applies to every maximal intersecting family
with covering number at least three. The structural result below
classifies equality under an additional shifting hypothesis.

\subsection{Shifting and equality cases}

We first define the shifting operation and the two properties used in the structural result.
Let $\mathcal{A} \subseteq \binom{[n]}{k}$ be a family. For $1 \leq i < j \leq n$ and $A \in \mathcal{A}$, define
\[
  s_{ij}(A)=
  \begin{cases}
  (A \setminus \{j\}) \cup \{i\}, & \text{if } j \in A,\ i \notin A,\ \text{and } (A \setminus \{j\}) \cup \{i\} \notin \mathcal{A};\\[2mm]
  A, & \text{otherwise.}
  \end{cases}
\]
Set $s_{ij}(\mathcal{A}) = \{s_{ij}(A) : A \in \mathcal{A}\}$. This operation is called the \emph{shifting (or shift) operation}, introduced in \cite{EKR}; see also \cite{frankl87}. We say that $\mathcal{A}$ is \emph{shifted} on $[n]$ if $s_{ij}(\mathcal{A}) = \mathcal{A}$ for all $1 \leq i < j \leq n$. More generally, $\mathcal{A}$ is shifted on a subset $P \subseteq [n]$ if $s_{ij}(\mathcal{A}) = \mathcal{A}$ for all $i < j$ with $i,j \in P$.

For a fixed integer $r>1$, a family $\mathcal A$ with
$\tau(\mathcal A)\ge r$ is called \emph{critically shifted} if it is
not shifted, but every shift that changes it reduces its covering
number below $r$. In this paper, critical shiftedness always refers
to the condition $\tau\ge3$. Both properties depend on the labeling
of the ground set; the theorem requires the existence of a labeling
with the stated property.

We now introduce two further types of construction for the equality
statement.
For $1\le r\le\lfloor k/2\rfloor$, choose a point $x$ and
pairwise disjoint sets
$A,B,C,D\subseteq[n]\setminus\{x\}$ such that
$|A|=k-r, \quad |B|=|C|=r, \quad |D|=k-2r$.
Put
$E_1=A\cup B, \quad E_2=A\cup C, \quad E_3=B\cup C\cup D$,
and define
\[
\mathcal R_r(n,k)
=
\{E_1,E_2,E_3\}
\cup
\left\{
F\in\binom{[n]}k:
x\in F,\;
F\cap E_i\ne\emptyset
\text{ for every }i\in[3]
\right\}.
\]
The isomorphism type of $\mathcal R_r(n,k)$ depends only on $n,k,r$.
Clearly, $\mathcal R_r(n,k)$ is maximal intersecting with covering
number three.
The three exceptional members have empty total intersection
and satisfy
$\bigl( |E_1\cap E_2|, |E_1\cap E_3|, |E_2\cap E_3| \bigr) = (k-r,r,r)$.
In particular,
$\mathcal R_1(n,k)\cong\mathcal G_1(n,k), \quad \mathcal R_2(n,k)\cong\mathcal Q^2(n,k)$.

Inclusion--exclusion gives
\begin{equation}\label{eq:R-size}
\begin{aligned}
|\mathcal R_r(n,k)|={}&\binom{n-1}{k-1}
-3\binom{n-k-1}{k-1}\\
&+\binom{n-k-r-1}{k-1}
+\binom{n-2k+r-1}{k-1}+3.
\end{aligned}
\end{equation}
Direct comparison with the formulas for $q$ and $g_1$ yields
\begin{equation}\label{eq:R-comparison}
|\mathcal R_r(n,k)|\le q(n,k)<g_1(n,k)
\qquad\left(2\le r\le\lfloor k/2\rfloor\right),
\end{equation}
where the first inequality is an equality exactly when $r=2$ or
$n=2k+1$. In particular, at $n=2k+1$, all $\mathcal R_r$ with
$r\ge2$ have size $q(2k+1,k)=\binom{2k}{k-1}-3k+3$.

For $k=4$, a further equality construction is needed.
Put $U=\{1,2,4\}$ and $L=\{3,5,6,7\}$, and define
\begin{equation}\label{eq:verified-special-four-description}
\mathcal F(n,4)=\{L\}\cup
\left\{F\in\binom{[n]}4:|F\cap U|\ge2,\ F\cap L\ne\emptyset\right\}.
\end{equation}
Clearly, this family is maximal intersecting with covering number
three. Counting according to $|F\cap U|$ gives
\[
|\mathcal F(n,4)|
=3\left(\binom{n-3}{2}-\binom{n-7}{2}\right)+5
=12n-61=q(n,4)=g_2(n,4).
\]

The constructions used in the equality statement admit critically
shifted labelings. For $\mathcal G_p$ use the labeling in its definition;
for $\mathcal R_r$, $r\ge2$, put the points in the order
$\{x\},A,B,C,D$, followed by the unused points. The family
$\mathcal F(n,4)$ uses its original labeling.
These assertions follow directly by checking which shifts change
the exceptional members.

\begin{thm}[Equality under critical shiftedness]\label{mre}
Let $n>2k$ and $k\ge4$.
Let $\mathcal F\subseteq\binom{[n]}k$ be maximal intersecting,
with $\tau(\mathcal F)\ge3$, and suppose that it admits a
critically shifted labeling.
Then either
$\mathcal F\cong\mathcal G_1(n,k)$ or
\[
|\mathcal F|\le M(n,k)
=
\max\{q(n,k),g_2(n,k)\}.
\]
Equality $|\mathcal F|=M(n,k)$ holds if and only if one of the
following occurs:
\begin{enumerate}
\item[\rm(i)]
$\mathcal F\cong\mathcal Q^2(n,k)$, where either $k=4$,
or $k\ge5$ and $n\le3k-2$;

\item[\rm(ii)]
$\mathcal F\cong\mathcal G_2(n,k)$, where either $k=4$,
or $k\ge5$ and $n\ge3k-2$;

\item[\rm(iii)]
$n=2k+1$ and
$\mathcal F\cong\mathcal R_r(2k+1,k)$ for some
$3\le r\le\lfloor k/2\rfloor$;

\item[\rm(iv)]
$k=4$ and either
$\mathcal F\cong\mathcal G_3(n,4)$ or
$\mathcal F\cong\mathcal F(n,4)$.
\end{enumerate}
\end{thm}

Distinct values of $r$ in part~{\rm(iii)} give nonisomorphic
families. Indeed, the hub $x$ is the unique point of maximum degree,
and an isomorphism must preserve the multiset $\{k-r,r,r\}$ of
pairwise intersection sizes of the exceptional members.
Thus the boundary $n=2k+1$ admits additional equality configurations
as $k$ increases; the theorem includes these alongside the exceptional
four-uniform constructions.

\subsection{Proof strategy and organization}

The local argument uses the exact sections determined by a minimum
three-cover and the common intersections of the single sections.
Replacing one single section by a suitable member and enlarging
the complementary pair section gives a controlled increase in size.
The cover criteria ensure that these replacements retain covering
number three. This reduces the extremal problem to the case of
three exceptional members and a small number of core configurations.

The second largest size requires an additional step because a shift
preserving $\tau\ge3$ need not preserve maximality. We prove a quantitative
inverse-shift bound using the connectivity of a disjointness graph.
It controls the size of a maximal family whose shifted image is
contained in $\mathcal G_1$, allowing the numerical bound to pass
from critically shifted families to all maximal families under
consideration. The equality classification is then obtained from
the same local replacements under the stated critical shiftedness
hypothesis.

Section~2 describes admissible shifting. Section~3 relates shifts
to minimum covers and introduces the exact sections. Section~4
collects the cross-intersection estimates, and Section~5 develops
the local replacements according to the core sizes. Section~6
proves the numerical maximum, the inverse-shift bound, and the two
main theorems.

Unless explicitly stated otherwise, assume $n>2k$ and $k\ge4$
in the structural arguments below.

\section{Admissible Shifting}

We now describe the shifting procedure used in the proof.
To bound the fully shifted case, we recall the following construction.
For $1\le s\le k$, define
\[
\mathcal K(n,k,s)
=
\left\{
A\in\binom{[n]}{k}:
1\in A \text{ and }
|A\cap[2,k+s-1]|\ge s-1
\right\}
\cup
\binom{[2,k+s-1]}{k}.
\]
The family $\mathcal K(n,k,s)$ is shifted and intersecting, with
$\tau(\mathcal K(n,k,s))=s$. Frankl and Wang proved that it gives
the sharp bound for shifted families with this covering constraint.

\begin{thm}[\cite{FW}]\label{fws}
Let $n, k, s$ be positive integers with $n > 2k \ge 2s$. If $\mathcal{F} \subseteq \binom{[n]}{k}$ is a shifted intersecting family with $\tau(\mathcal{F}) \ge s$, then
\[
|\mathcal{F}| \le |\mathcal{K}(n,k,s)|.
\]
\end{thm}

For $s=3$, counting gives
\[
|\mathcal K(n,k,3)|
=
\binom{n-1}{k-1}
-
k\binom{n-k-2}{k-2}
-
\binom{n-k-1}{k-1}
+
(k+1).
\]
Direct comparison with $g_2$ yields
\begin{align}\label{gk3}
|\mathcal G_2(n,k)|
-
|\mathcal K(n,k,3)|
=
(k-3)\binom{n-k-2}{k-2}
-
\sum_{i=k+4}^{2k}\binom{n-i}{k-2}
-
(k-3)
>0.
\end{align}

Call a shift of a family $\mathcal F$ with $\tau(\mathcal F)\ge3$
\emph{admissible} if $\tau(s_{ij}(\mathcal F))\ge3$.
Such a shift preserves intersection and cardinality. Starting with
a maximum intersecting family subject to $\tau\ge3$, apply changing
admissible shifts until none remains. Each changing shift strictly
decreases $\sum_{F\in\mathcal F}\sum_{x\in F}x$, so the process
terminates, and every intermediate family remains maximum and
hence maximal.

Since the maximum size is at least $g_1>g_2$, Theorem~\ref{fws}
and~\eqref{gk3} exclude a fully shifted terminal family.
The terminal family is therefore critically shifted. This supplies
the reduction used for the numerical maximum.

For a merely maximal family, admissible shifting need not preserve
maximality. In Section~6, the inverse-shift bound controls maximal
completions of shifted images and justifies the corresponding
reduction for the second largest size.

\section{Reduction to Cross-Intersecting Families}

We relate the obstruction to an admissible shift to minimum covers,
and then use a three-cover to define the exact sections.
For a family $\mathcal{F} \subseteq \binom{[n]}{k}$ with $\tau(\mathcal{F}) = r$, define
\[
\mathcal{T}(\mathcal{F}) = \left\{ X \in \binom{[n]}{r} : X \text{ is a cover of } \mathcal{F} \right\}.
\]

The following proposition describes a shift that lowers the covering number.

\begin{prop}
Let $\mathcal{F} \subseteq \binom{[n]}{k}$ be an intersecting family with $\tau(\mathcal{F}) \ge 3$. If $\tau(s_{ij}(\mathcal{F})) < 3$ for some $i,j \in [n]$, then
\[
\tau(\mathcal{F}) = 3, \qquad \tau(s_{ij}(\mathcal{F})) = 2,
\]
and there exists $a \in [n]$ such that $\{a,i,j\} \in \mathcal{T}(\mathcal{F})$.
\end{prop}
\begin{proof}
Without loss of generality, assume $\tau(s_{12}(\mathcal{F})) < 3$, and let $T$ be a minimum cover of $s_{12}(\mathcal{F})$. Then $|T| \le 2$.

If $|T| = 1$, say $T = \{x\}$, then $\{x, 2\}$ would be a cover of $\mathcal{F}$, contradicting $\tau(\mathcal{F}) \ge 3$. Hence $|T| = 2$.
Moreover, we must have $1 \in T$. Otherwise, $T$ would also be a cover of $\mathcal{F}$, which is again a contradiction. Thus, $T = \{1, a\}$ for some $a \in [n]$.
Finally, for any $F \in \mathcal{F}$, either $F \in s_{12}(\mathcal{F})$, or $1 \notin F$ and $2 \in F$, in which case $F \cap \{1, 2, a\} \neq \emptyset$. Consequently, $\{1, 2, a\}$ is a cover of $\mathcal{F}$, which implies $\tau(\mathcal{F}) = 3$ and $\tau(s_{12}(\mathcal{F})) = 2$.
\end{proof}

This demonstrates that every obstruction to shifting is witnessed by a 3-element transversal. As an immediate consequence, we obtain the following result.

\begin{prop}\label{p2}
Let $\mathcal{F} \subseteq \binom{[n]}{k}$ be critically shifted
with $\tau(\mathcal{F})\ge3$. Then $\tau(\mathcal F)=3$.
If no set $T \in \mathcal{T}(\mathcal{F})$ contains $\{i,j\}$,
then $s_{ij}(\mathcal{F}) = \mathcal{F}$.
\end{prop}

Thus a changing admissible shift is impossible in a critically shifted family.

\medskip

\begin{lem}\label{lem:verified-cover-labeling}
A critically shifted family with covering number at least three
admits a critically shifted labeling for which $[3]$ is a cover.
\end{lem}

\begin{proof}
Choose a three-cover $T$ minimizing $\sum_{t\in T}t$.
If $i<j$, $i\notin T$ and $j\in T$, then $s_{ij}(\mathcal F)$
changes the family: otherwise $(T\setminus\{j\})\cup\{i\}$
would be a cover with smaller sum. Hence every such shift reduces
the covering number to two.

Place the points of $T$ first, preserving the relative order within
$T$ and within its complement. The only reversed pairs are those
just considered. For the transposition $\pi_{ij}$,
$s_{ji}(\mathcal F)=\pi_{ij}(s_{ij}(\mathcal F))$, so reversing
these pairs preserves the covering number of their shifted images.
All other pair orientations are unchanged. Thus the new labeling
is critically shifted and has cover $[3]$.
\end{proof}

Let $\mathcal F$ be maximal intersecting with minimum cover $[3]$.
For $P\subseteq[3]$, define
\[
\mathcal F\langle P\rangle=\{F\in\mathcal F:F\cap[3]=P\},
\qquad
\mathcal F(P)=\{F\setminus P:F\in\mathcal F\langle P\rangle\}.
\]
We suppress braces and commas, writing $\mathcal F(12)$ and
$\mathcal F(123)$, for example. Then
\[
\mathcal F=\bigcup_{\emptyset\ne P\subseteq[3]}\mathcal F\langle P\rangle,
\qquad \mathcal F(123)=\binom{[4,n]}{k-3}.
\]
Each $\mathcal F(i)$ is nonempty, since $[3]$ is a minimum cover.
For $\{i,j,\ell\}=[3]$, maximality gives
\begin{equation}\label{eq:completion-sections}
\mathcal F(j\ell)=\left\{B\in\binom{[4,n]}{k-2}:
B\cap A\ne\emptyset\text{ for every }A\in\mathcal F(i)\right\}.
\end{equation}
The reverse saturation of $\mathcal F(i)$ also involves
$\mathcal F(j)$ and $\mathcal F(\ell)$; it does not follow from
$\mathcal F(j\ell)$ alone.
Put
\[
X_i=\bigcap_{A\in\mathcal F(i)}A,
\qquad Y_i=\bigcap_{B\in\mathcal F([3]\setminus\{i\})}B.
\]
Here the intersection of an empty subfamily of $2^{[4,n]}$ is
interpreted as $[4,n]$. In particular, $Y_i$ is defined even when
the corresponding pair section is empty.
For a family $\mathcal C$ and a point $z$ outside the fixed cover,
write
\[
\mathcal C(z)=\{C\setminus\{z\}:z\in C\in\mathcal C\},
\qquad \mathcal C(\overline z)=\{C\in\mathcal C:z\notin C\}.
\]
Thus $\mathcal F(i,z)$ means $(\mathcal F(i))(z)$, and similarly
for $\mathcal F(i,\overline z)$.
We also write $S+z=S\cup\{z\}$ and
$\mathcal F(\overline x)=\{F\in\mathcal F:x\notin F\}$ when the
entire family, rather than one exact section, is being restricted.

By \eqref{eq:completion-sections}, every $(k-2)$-set meeting $X_i$
belongs to the complementary pair section. Consequently,
$|X_i|\ge2$ implies $Y_i=\emptyset$, while $|X_i|=1$ implies
$|Y_i|\le1$. A pair $\{i,z\}$ covers $\mathcal F$ precisely when
$z\in X_j\cap X_\ell\cap Y_i$.
These observations will be used to check the covering number after
a local replacement.

\section{Some Results on Cross-intersecting Families}

We collect the cross-intersection bounds needed for the exact sections.
Neighborhoods and shadows connect these estimates with the
Kruskal--Katona theorem.

\subsection{Neighborhoods and Shadows}
For a family $\mathcal{A} \subseteq \binom{[n]}{k}$, we define its $\ell$-neighborhood and $\ell$-shadow by
\[
N_{\ell}^{[n]}(\mathcal{A})
=
\Bigl\{\,B\in \binom{[n]}{\ell}:\ B\cap A=\emptyset
\text{ for some } A\in \mathcal{A}\,\Bigr\},
\]
and
\[
\partial_{\ell}(\mathcal{A})
=
\Bigl\{\,B\in \binom{[n]}{\ell}:\ B\subset A
\text{ for some } A\in \mathcal{A}\,\Bigr\}.
\]
When the ground set $[n]$ is clear from the context, we simply write $N_{\ell}(\mathcal{A})$.

Let
\[
\overline{\mathcal{A}}
=
\{\, [n]\setminus A:\ A\in\mathcal{A}\,\}.
\]
Then we have the duality identity
\[
N_{\ell}(\mathcal{A})=\partial_{\ell}(\overline{\mathcal{A}}),
\]
which allows us to convert neighborhood estimates into shadow estimates.

\subsection{Kruskal--Katona Theorem}

We next recall the classical Kruskal--Katona theorem, which provides sharp lower bounds on shadows. To state it, we first introduce the $k$-binomial expansion (also called the $k$-cascade). For positive integers $m$ and $k$, there exist unique integers
\[
a_1>a_2>\cdots>a_t\ge k-t+1\ge 1
\]
such that
\[
m
=
\binom{a_1}{k}
+\binom{a_2}{k-1}
+\cdots
+\binom{a_t}{k-t+1}.
\]
We refer to $(a_1,\ldots,a_t)$ as the $k$-cascade of $m$.

To describe extremal configurations, we introduce two natural orderings on $\binom{[n]}{k}$. For distinct sets $F,G\in\binom{[n]}{k}$, we write
\[
F\prec_S G \quad\text{if}\quad \max(F\setminus G)<\max(G\setminus F),
\]
and
\[
F\prec_L G \quad\text{if}\quad \min(F\setminus G)<\min(G\setminus F).
\]
Let $\mathcal{S}(k,m)$ denote the first $m$ $k$-sets in the squashed order, and let $\mathcal{L}([n],k,m)$ denote the initial segment of size $m$ in $\binom{[n]}{k}$ with respect to $\prec_L$.

The Kruskal--Katona theorem can now be stated as follows.

\begin{thm}[Kruskal \cite{Kruskal} and Katona \cite{Katona68}]
\label{thm:KK}
Let $\mathcal{A}\subseteq \binom{[n]}{k}$, and write $|\mathcal{A}|$ in its $k$-binomial expansion as
\[
|\mathcal{A}|
=
\binom{a_1}{k}
+\binom{a_2}{k-1}
+\cdots
+\binom{a_t}{k-t+1}.
\]
Then, for every $1 \le \ell \le k$,
\[
\bigl|\partial_{\ell}(\mathcal{A})\bigr|
\ge
\binom{a_1}{\ell}
+\binom{a_2}{\ell-1}
+\cdots
+\binom{a_t}{\ell-t+1},
\]
where terms with negative lower index are interpreted as zero. Moreover, equality holds if $\mathcal{A}\cong\mathcal{S}(k,|\mathcal{A}|)$.
\end{thm}

The dual form used below is the following.
\begin{thm}[co-Kruskal--Katona]\label{lem1}
Let $n\ge k+\ell$ and let $\mathcal A\subseteq\binom{[n]}k$.
Then
\[
|N_\ell(\mathcal A)|\ge
|N_\ell(\mathcal L([n],k,|\mathcal A|))|.
\]
\end{thm}
We use the uniqueness criterion of F\"uredi--Griggs and M\"ors
only in the specific instance verified in the proof of
Lemma~\ref{zh}; equality in the neighborhood bound alone does
not imply uniqueness in general.

\subsection{Classical bounds for cross-intersecting families}

We begin with the following theorem of Frankl and Tokushige.

\begin{thm}[Frankl--Tokushige \cite{FT}]\label{ft}
Let $n,k,\ell$ be positive integers with $n\ge k+\ell$ and $\ell\ge k$. Let $\mathcal A\subseteq\binom{[n]}k$ and $\mathcal B\subseteq\binom{[n]}\ell$ be nonempty cross-intersecting families. Then
\[
|\mathcal A|+|\mathcal B|
\le\binom n\ell-\binom{n-k}\ell+1.
\]
Unless $n=k+\ell$ or $k=\ell=2$, the inequality is strict whenever $|\mathcal A|>1$ and $|\mathcal B|>1$.
\end{thm}

For nontrivial cross-intersecting families, we shall use the following bound.

\begin{thm}[Frankl \cite{fr24}]\label{cab1}
Let $n,k,\ell$ be positive integers with $n\ge k+\ell$ and $\ell\ge k\ge2$. Let $\mathcal A\subseteq\binom{[n]}k$ and $\mathcal B\subseteq\binom{[n]}\ell$ be nonempty, nontrivial cross-intersecting families. Then
\[
|\mathcal A|+|\mathcal B|
\le\binom n\ell-2\binom{n-k}\ell
+\binom{n-2k}\ell+2.
\]
\end{thm}

For unequal uniformities, the following result of Frankl and Wang
requires only the family of smaller uniformity to be nontrivial.

\begin{thm}[Frankl--Wang \cite{FW}]\label{cab2}
Let $n>k+\ell$ and $\ell>k\ge2$.
Let $\mathcal A\subseteq\binom{[n]}k$ and
$\mathcal B\subseteq\binom{[n]}\ell$ be nonempty cross-intersecting
families, with $\bigcap_{A\in\mathcal A}A=\emptyset$.
Then
\[
|\mathcal A|+|\mathcal B|
\le\binom n\ell-2\binom{n-k}\ell+\binom{n-2k}\ell+2.
\]
Equality holds exactly when $\mathcal A$ consists of two disjoint
$k$-sets and $\mathcal B$ consists of all $\ell$-sets meeting both.
\end{thm}

We also use the following consequence of M\"ors's shadow theorem;
see~\cite{mh} and~\cite[Section 3]{FT98}.

\begin{thm}[M\"ors~\cite{mh}]\label{mh}
Let $k,\ell\ge2$, $n>k+\ell$, and $n\ge2\max\{k,\ell\}$.
Suppose that $\mathcal A\subseteq\binom{[n]}k$ and
$\mathcal B\subseteq\binom{[n]}\ell$ are nonempty, nontrivial
cross-intersecting families. Put
\[
\alpha=\binom{n-1}{k-1}-\binom{n-\ell-1}{k-1}+1,
\qquad
\beta=\binom{n-1}{\ell-1}-\binom{n-k-1}{\ell-1}+1.
\]
Then $|\mathcal A|\le\alpha$ or $|\mathcal B|\le\beta$.
If either inequality fails, the other holds strictly.
\end{thm}

\subsection{Further cross-intersection estimates}

Removing the common intersection of one family yields the following
refinement.

\begin{thm}[A bound in terms of the common intersection]\label{zh1}
Let $n,k,\ell$ be positive integers with $n>k+\ell$ and $\ell\ge3$. Let $\mathcal A\subseteq\binom{[n]}k$ and $\mathcal B\subseteq\binom{[n]}\ell$ be nonempty cross-intersecting families. Put
\[
C=\bigcap_{A\in\mathcal A}A,\qquad a=|C|,
\qquad
\mathcal B_0=\{B\in\mathcal B:B\cap C=\emptyset\}.
\]
If $k-\ell<a<k$ and $\mathcal B_0\ne\emptyset$, then
\[
|\mathcal A|+|\mathcal B|
\le\binom n\ell-2\binom{n-k}\ell
+\binom{n-2k+a}\ell+2.
\]
Equality holds if and only if $\mathcal A=\{A_1,A_2\}$ with $|A_1\cap A_2|=a$ and
\[
\mathcal B=\left\{B\in\binom{[n]}\ell:
B\cap A_i\ne\emptyset\text{ for each }i\in[2]\right\}.
\]
\end{thm}

\begin{proof}
Put $\mathcal A'=\{A\setminus C:A\in\mathcal A\}$.
Then $\mathcal A'$ is nontrivial and cross-intersecting with
$\mathcal B_0$, and $1\le k-a<\ell$. Moreover,
\begin{equation}\label{eq:zh1-decomposition}
|\mathcal B|\le\binom n\ell-\binom{n-a}\ell+|\mathcal B_0|,
\end{equation}
with equality precisely when $\mathcal B$ contains every $\ell$-set
meeting $C$. If $k-a\ge2$, Theorem~\ref{cab2} gives
\begin{equation}\label{eq:zh1-reduced}
|\mathcal A'|+|\mathcal B_0|
\le\binom{n-a}\ell-2\binom{n-k}\ell
+\binom{n-2k+a}\ell+2.
\end{equation}
Equality holds precisely when $\mathcal A'$ consists of two disjoint
sets and $\mathcal B_0$ is their full cross-family on $[n]\setminus C$.

If $k-a=1$, write $\mathcal A'=\{\{x_1\},\ldots,\{x_t\}\}$.
Since $\mathcal B_0\ne\emptyset$, we have $2\le t\le\ell$, and
\[
|\mathcal A'|+|\mathcal B_0|
\le t+\binom{n-a-t}{\ell-t}
\le2+\binom{n-a-2}{\ell-2}.
\]
The second inequality is strict for $t>2$, since $n-a\ge\ell+2$.
Pascal's identity now gives~\eqref{eq:zh1-reduced} with the same
equality condition.

Combining~\eqref{eq:zh1-decomposition} and~\eqref{eq:zh1-reduced}
proves the bound. Their equality conditions give the stated
configuration, whose size follows directly by inclusion--exclusion.
\end{proof}

\subsection{Small cross-intersecting families}

We begin with graph estimates for the smallest uniformities, including
the six-point ground set arising when $n=9$ and $k=4$.

\begin{lem}\label{lem:completion-graphs}
Let $N\ge6$.
\begin{enumerate}
\item[\rm(i)] Nonempty cross-intersecting graphs on $M\ge4$ vertices
have total size at most $2M-2$. If both are nontrivial, their total
size is at most six.
\item[\rm(ii)] Suppose $\mathcal A\subseteq\binom{[N]}3$ is nontrivial,
$|\mathcal A|\le3N-8$, and
$\mathcal B\subseteq\binom{[N]}2$ is cross-intersecting with it.
Then $|\mathcal A|+|\mathcal B|\le3N-5$, with equality only if
$|\mathcal A|=3N-8$. If $|\mathcal B|\ge4$, the stronger bound
$|\mathcal A|+|\mathcal B|\le2N$ holds.
\end{enumerate}
\end{lem}
\begin{proof}
Part (i) follows directly from Theorems~\ref{ft} and~\ref{cab1}
with both uniformities equal to two.

For (ii), view $\mathcal B$ as a graph $G$ and write $m=e(G)$.
Every member of $\mathcal A$ is a three-vertex cover of $G$.
Choosing such a cover avoiding each prescribed vertex gives
$\Delta(G)\le3$, and hence $m\le9$.
For $m\le3$ the assertion is immediate, so suppose $m\ge4$.

If $G$ has a matching of size three, each three-cover chooses one
endpoint of each matching edge. There are eight choices.
No edge can join a matching endpoint to a vertex outside this
matching, since it would put that endpoint in every member of
$\mathcal A$; an edge entirely outside it is also impossible.
Put $r=m-3\le6$. Each of the $r$ additional edges excludes the two
endpoints of a different edge of the three-dimensional cube of
these eight choices. Any $r\le6$ cube edges have at least $r$
endpoints: a subgraph on at most five cube vertices has at most
as many edges as vertices. Thus $|\mathcal A|\le8-r$, and
$|\mathcal A|+m\le11\le2N$.

Otherwise the matching number is two. Suppose first that the
vertex-cover number is three. If $G$ is triangle-free, it is not
bipartite, by K\"onig's theorem. It therefore contains a five-cycle;
a longer odd cycle or any further edge incident with or outside
this cycle would give a matching of size three, unless a chord
created a triangle. Hence $G$ is a five-cycle with isolated vertices,
and there are five three-covers.
If $G$ contains a triangle, every triangle vertex has at most one
neighbor outside it. Write $b\le3$ for the number of edges leaving
the triangle. If there is no edge outside the triangle, then
$b\ge1$ and the number of three-covers is
$1+(3-b)(N-3)+b$. Adding $3+b$ gives at most $2N$.
Otherwise the graph outside the triangle has matching number one.
It cannot contain a triangle, since that would require four vertices
to cover $G$, so it is a star of at most three edges.
If it has one edge, the number of three-covers is at most $6-b$,
giving a total at most ten. If it has two or three edges, every
three-cover consists of its center and two triangle vertices.
There are at most three choices, so the total is at most twelve.

Finally, suppose $\{u,v\}$ is a two-cover of $G$. Put
$e=\mathbf1_{uv\in G}$, $a=|N_G(u)\setminus\{v\}|$,
$b=|N_G(v)\setminus\{u\}|$, and
$q=|(N_G(u)\cup N_G(v))\setminus\{u,v\}|$.
Here $a,b\ge1$, $a,b\le3-e$, and $m=a+b+e$.
The number of three-covers is
\[
N-2+f(a)+f(b)+(1-e)\binom{N-2-q}{3-q},
\qquad f(1)=N-3,\quad f(2)=1,\quad f(3)=0.
\]
For $m=4,5,6$ this is at most $2N-4,N,N-1$, respectively,
as follows by substituting $1\le a,b\le3-e$ and
$\max\{a,b\}\le q\le a+b$ (with $q\ge2$ if $e=0$).
Adding $m$ gives at most $2N$. This proves (ii).
\end{proof}

The next estimate concerns families whose uniformities differ by one,
with a size restriction on the family of larger sets.

\begin{lem}\label{lem:completion-small-section}\label{zh}
Let $t\ge2$, $v>2t+1$, and let
$\mathcal A\subseteq\binom{[v]}{t+1}$ be nontrivial and
$\mathcal B\subseteq\binom{[v]}t$ be cross-intersecting with $\mathcal A$.
If $|\mathcal A|\le H=\binom{v-1}t-\binom{v-t-2}t+1$, then
\[
|\mathcal A|+|\mathcal B|\le\binom vt-\binom{v-t-1}t+1.
\]
Equality requires $|\mathcal A|=H$.
\end{lem}
\begin{proof}
For $t=2$, the assertion is Lemma~\ref{lem:completion-graphs}(ii)
with $N=v$. Hence assume $t\ge3$.
Put $\beta=\binom{v-1}{t-1}-\binom{v-t-2}{t-1}+1$ and
$\alpha=\binom{v-1}t-\binom{v-t-1}t+1$.
If $|\mathcal B|<\beta$, combining this inequality with
$|\mathcal A|\le H$ proves the bound and its equality condition. Otherwise $\mathcal B$ is nontrivial. Indeed, if every member of
$\mathcal B$ contained a point $x$, a member of $\mathcal A$ avoiding
$x$ would give $|\mathcal B|\le\beta-1$.
M\"ors's theorem gives $|\mathcal A|\le\alpha$; if equality holds,
then $|\mathcal B|\le\beta$, and the desired inequality is strict.

Suppose $|\mathcal A|\le\alpha-1$ and replace it by its lexicographic
initial segment, replacing $\mathcal B$ by the full cross-family.
All members of the first family contain $1$. The members of the
second family avoiding $1$ form a nonempty $t$-uniform family.
The Frankl--Tokushige inequality, after deleting $1$ from the first
family, gives the desired bound. Equality would require the second
family avoiding $1$ to be a singleton and $|\mathcal A|=\alpha-1$.
Equality must also hold in the corresponding neighborhood bound.

Here the complements of $\mathcal A$ are $p$-sets, where $p=v-t-1>t$,
and their $p$-binomial expansion is
\[
\alpha-1=\sum_{j=0}^{t-1}\binom{v-2-j}{p-j}.
\]
Equality in their $t$-shadow forces equality in their immediate shadow.
If $p=t+1$, these are the same shadow, so this is automatic.
If $p\ge t+2$, the minimum immediate shadow has the
$(p-1)$-binomial expansion
$\sum_{j=0}^{t-1}\binom{v-2-j}{p-1-j}$.
Increasing this minimum by one appends the valid term
$\binom{p-t-1}{p-t-1}$.
Applying Kruskal--Katona for another $p-t-1$ shadow operations
then increases the lower bound for the $t$-shadow by
$\binom{p-t-1}{0}=1$, a contradiction.
The uniqueness criterion of F\"uredi and Griggs~\cite{FG}
and M\"ors~\cite{mh}, as stated in~\cite[Theorem 2]{SV}, applies
because the expansion has $t<p$ terms.
Thus $\mathcal A$ would be isomorphic to the lexicographic family,
which is trivial, a contradiction.
\end{proof}

We also need the following estimate when the uniformities differ by two.

\begin{lem}\label{lem:completion-small-family}
Let $s\ge3$, $v\ge2s+1$, and let nonempty cross-intersecting families
$\mathcal A\subseteq\binom{[v]}{s+1}$ and
$\mathcal B\subseteq\binom{[v]}{s-1}$ satisfy $|\mathcal A|\le v-s$.
Then $|\mathcal A|+|\mathcal B|\le
\binom v{s-1}-\binom{v-s-1}{s-1}+1$.
\end{lem}
\begin{proof}
Write $h=v-s$ and $u=|\mathcal A|$.
Its lexicographic replacement consists of a fixed $s$-set together
with $u$ different additional points. The full cross-family has size
$\binom v{s-1}-\binom h{s-1}+\binom{h-u}{s-1-u}$.
For $1\le u\le s-1$, the function
$u+\binom{h-u}{s-1-u}$ is nonincreasing; for $s\le u\le h$ it is
$u\le h\le1+\binom{h-1}{s-2}$.
Kruskal--Katona and Pascal's identity prove the result.
\end{proof}

\section{The Core Sets $X_i$}

Throughout this section $n>2k$, $k\ge4$, and $\mathcal F$ is maximal intersecting with
minimum cover $[3]$. Critical shiftedness and a size assumption
are imposed where needed. We estimate the complementary section
sums and use these estimates to choose local replacements that
preserve covering number three.

In this section, critical shiftedness is used only through the
following consequence of Proposition~\ref{p2}: if $i\in[3]$ and
$z\in[4,n]$ belong to no common three-cover, then
$s_{iz}(\mathcal F)=\mathcal F$. This property is preserved by
any permutation of the three cover points, since such a permutation
preserves the direction of every shift from $[4,n]$ into $[3]$.
Thus the arguments below remain valid after the indicated permutations
of $[3]$, even if the resulting labeling is no longer critically shifted.

For the repeated size comparisons, put
\[
K=\binom{n-3}{k-2}-\binom{n-k-2}{k-2}+1,
\qquad \delta=\binom{n-k-3}{k-3}-1=g_1(n,k)-g_2(n,k)>0.
\]
This notation is retained in Section~6. When $n,k$ are fixed, we
abbreviate $g_p(n,k)$ and $q(n,k)$ to $g_p$ and $q$.

\subsection{Extremal threshold inequalities}
For the numerical comparisons, extend the construction
$\mathcal F(n,4)$ to general $k$ as follows.
Define $\mathcal F(n,k)=\mathcal M\cup\mathcal N$, where
\[
\begin{aligned}
\mathcal M
={}&\left\{B\in\binom{[2,n]}{k}: \{2,4\}\subset B,
  \ 3\in B\text{ or }B\cap[5,k+3]\ne\emptyset\right\}\\
&\quad\cup\{[3,k+3]\setminus\{4\}\},
\end{aligned}
\]
\[
\mathcal N
=
\left\{A\in\binom{[n]}{k}:
1\in A \text{ and }
A\cap B\neq\emptyset \text{ for all } B\in\mathcal M
\right\}.
\]

Counting its sections gives
\[
|\mathcal F(n,k)|
=
\binom{n-3}{k-3}
+
3\left(
\binom{n-3}{k-2}
-
\binom{n-k-2}{k-2}
+
1
\right)
+
2\left(
\binom{n-k-3}{k-3}
-
1
\right).
\]
Direct comparison of this expression and the formula for $g_2$
gives the following estimates.
\begin{lem}\label{eq12}
Let $n>2k$ and $k\ge4$. Then
\[
\binom{n-3}{k-3}
+3\Bigl(\binom{n-3}{k-2}-\binom{n-k-2}{k-2}+1\Bigr)
+2\Bigl(\binom{n-k-3}{k-3}-1\Bigr)
\le |\mathcal G_2(n,k)|.
\]
Equality holds if and only if $k=4$.
\end{lem}

\begin{proof}
Let $D(n,k)$ be the difference between the right and left sides.
Substituting the formula for $g_2(n,k)$ and using Pascal's identity gives
\begin{align*}
D(n,k)
={}&\binom{n-3}{k-1}-\binom{n-3}{k-2}
-2\binom{n-k-2}{k-1}\\
&-3\binom{n-k-3}{k-3}+\binom{n-2k}{k-1}+3.
\end{align*}
For $k=4$, this expression is zero for every $n\ge9$.
Suppose henceforth that $k\ge5$. At $n=2k+1$, we have
\[
D(2k+1,k)=\frac1k\binom{2k-2}{k-1}-3k+7>0.
\]
Indeed, the binomial term equals $14$ when $k=5$ and is multiplied
by $(4k-2)/(k+1)\ge3$ when $k$ is replaced by $k+1$;
thus it exceeds $3k-7$ for every $k\ge5$.

For $n\ge2k+2$, another application of Pascal's identity yields
\begin{align*}
D(n,k)-D(n-1,k)-D(n-1,k-1)
={}&2\binom{n-k-3}{k-3}+3\binom{n-k-4}{k-5}\\
&-\binom{n-2k}{k-3}-\binom{n-2k-1}{k-3}-3
\ge0.
\end{align*}
The last inequality follows from $n-k-3\ge n-2k$
and $\binom{n-k-4}{k-5}\ge1$.
Induction on $n$, using the cases $k=4$ and $n=2k+1$ above,
now gives $D(n,k)>0$ whenever $k\ge5$ and $n>2k$.
Therefore the required inequality holds, with equality exactly when
$k=4$.
\end{proof}

\begin{lem}\label{lem:completion-strong-number}
For $k\ge5$ and $n>2k$,
\[
\binom{n-3}{k-3}+3(K+\delta)<g_2(n,k).
\]
\end{lem}
\begin{proof}
	Let
	\[
	D(n,k)=g_2(n,k)-\binom{n-3}{k-3}-3(K+\delta).
	\]
	Substituting the formula for $g_2(n,k)$ and using Pascal's identity gives
	\[
	D(n,k)
	=\binom{n-3}{k-1}-\binom{n-3}{k-2}
	-2\binom{n-k-2}{k-1}
	-4\binom{n-k-3}{k-3}
	+\binom{n-2k}{k-1}+4.
	\]
	For $k=5$, this simplifies to
	\[
	D(n,5)=2(n-10)(n-9)>0
	\qquad(n\ge11).
	\]
	At $n=2k+1$, we have
	\[
	D(2k+1,k)=\frac1k\binom{2k-2}{k-1}-4k+10>0
	\qquad(k\ge5).
	\]
	Indeed, the first term equals $14$ when $k=5$ and is multiplied by
	$(4k-2)/(k+1)\ge3$ when $k$ is replaced by $k+1$; hence it exceeds
	$4k-10$ for every $k\ge5$.
	
	For $k\ge6$ and $n\ge2k+2$, another application of Pascal's identity yields
	\[
	\begin{aligned}
		&D(n,k)-D(n-1,k)-D(n-1,k-1)\\
		&\qquad=
		2\binom{n-k-3}{k-3}
		+4\binom{n-k-4}{k-5}
		-\binom{n-2k}{k-3}
		-\binom{n-2k-1}{k-3}-4
		\ge0.
	\end{aligned}
	\]
	The last inequality follows from $n-k-3\ge n-2k$ and
	$\binom{n-k-4}{k-5}\ge1$.
	Induction on $n$, using the cases $k=5$ and $n=2k+1$ above, now gives
	$D(n,k)>0$ whenever $k\ge5$ and $n>2k$, as required.
\end{proof}

\subsection{Decomposition estimates}
The case $k=4$ requires a separate argument when a core has size two.
\begin{lem}\label{lem:maximum-k4-core-two}
Let $n\ge9$ and let $\mathcal F\subseteq\binom{[n]}4$ be maximal
intersecting, with $\tau(\mathcal F)=3$ and $|\mathcal F|\ge12n-61$.
Suppose that $[3]$ is a minimum cover and
\[
\left|\bigcap_{E\in\mathcal F(1)}E\right|=2.
\]
Then $|\mathcal F(1)|\in\{2,3\}$, and
\[
|\mathcal F(1)|+|\mathcal F(23)|
=
2n-6
=
\binom{n-3}{2}
-\binom{n-6}{2}
-\binom{n-7}{1}+2.
\]
\end{lem}

\begin{proof}
By relabeling $[4,n]$, we may assume that $X_1=\{4,5\}$.
Put $m=|\mathcal F(1)|\ge2$.
Maximality gives
\begin{equation}\label{eq:max-k4-full-cross}
\mathcal F(23)
=
\left\{
E\in\binom{[4,n]}2:
E\cap A\ne\emptyset
\text{ for every }A\in\mathcal F(1)
\right\}.
\end{equation}

Suppose, for a contradiction, that $m\ge4$.
Every member of $\mathcal F(2)\cup\mathcal F(3)$ meets $\{4,5\}$,
since a three-set avoiding this pair cannot meet all $m$ members
of $\mathcal F(1)$.
Equation \eqref{eq:max-k4-full-cross} also gives
\[
\mathcal F(23)
=
\left\{
E\in\binom{[4,n]}2:E\cap\{4,5\}\ne\emptyset
\right\}.
\]
Thus every member of $\mathcal F$ avoiding $1$ meets
$\{4,5\}$.
By maximality, every four-set containing $\{1,4,5\}$
belongs to $\mathcal F$.
In particular,
\[
\mathcal F(1)
=
\left\{
E\in\binom{[4,n]}3:\{4,5\}\subseteq E
\right\}.
\]

It follows that $\mathcal F$ contains every four-set
containing one of the triples
\[
\{1,2,3\},\qquad
\{1,4,5\},\qquad
\{2,3,4\},\qquad
\{2,3,5\}.
\]
Since $n\ge9$ and $\mathcal F$ is intersecting, every member
of $\mathcal F$ meets all four triples.
Consequently, every member contains at least two points
of $[5]$, and its intersection with $[5]$, when of size
two, belongs to
\[
\mathcal E
=
\bigl\{
\{1,2\},\{1,3\},
\{2,4\},\{2,5\},
\{3,4\},\{3,5\}
\bigr\}.
\]

For each $e\in\mathcal E$, define
\[
\mathcal H_e
=
\left\{
F\setminus[5]:
F\in\mathcal F,\ F\cap[5]=e
\right\}
\subseteq\binom{[6,n]}2
\]
and
\[
p_e
=
|\mathcal H_e|
+
\left|
\left\{
x\in[6,n]:
([5]\setminus e)\cup\{x\}\in\mathcal F
\right\}
\right|.
\]
Maximality identifies the set counted by the second term with
$\bigcap_{H\in\mathcal H_e}H$ if $\mathcal H_e\ne\emptyset$,
and with $[6,n]$ otherwise.

Every four-subset of $[5]$ belongs to $\mathcal F$.
Moreover, the four triples displayed above account for
all three-subsets of $[5]$ whose complements do not
belong to $\mathcal E$.
Therefore,
\begin{equation}\label{eq:max-k4-five-point-count}
|\mathcal F|
=
4(n-5)+5+\sum_{e\in\mathcal E}p_e.
\end{equation}

We claim that if $p_e>n-5$, then
\begin{equation}\label{eq:max-k4-large-section}
p_e\le2n-13,
\end{equation}
and some $f\in\mathcal E$ disjoint from $e$ satisfies
$p_f=3$.

The assumption $p_e>n-5$ forces $\mathcal H_e$ to be nonempty
with empty total intersection, so $p_e=|\mathcal H_e|>n-5$.
Since $\tau(\mathcal F)=3$, the pair $e$ does not cover
$\mathcal F$.
The empty total intersection of $\mathcal H_e$ means that
no member of $\mathcal F$ has intersection $[5]\setminus e$
with $[5]$.
Hence there is some $f\in\mathcal E$ disjoint from $e$
such that $\mathcal H_f\ne\emptyset$.

The families $\mathcal H_e$ and $\mathcal H_f$ are
cross-intersecting.
Two distinct pairs on $[6,n]$ have at most $n-5$ common
intersecting pairs: the number is $n-5$ when they intersect
and $4\le n-5$ when they are disjoint.
Since $|\mathcal H_e|>n-5$, the family $\mathcal H_f$
must therefore be a singleton.
Consequently, $p_f=3$, and every member of $\mathcal H_e$
meets its unique pair.
This gives $p_e\le2n-13$, proving the claim.

In particular, the edges $e\in\mathcal E$ satisfying
$p_e>n-5$ are pairwise intersecting.
Since $\mathcal E$ is the edge set of $K_{2,3}$, there
are at most three such edges.
Put $d=n-5\ge4$, and let $L$ be the number of edges with
$p_e>d$. Each such edge satisfies $p_e\le2d-3$ and has a
disjoint witness $f$ with $p_f=3$.
If $L=0$, the sum of the six terms is at most $6d$.
If $L=1$, it is at most $(2d-3)+3+4d=6d$.
If $L=2$, at least one other term equals $3$, giving
\[
\sum_{e\in\mathcal E}p_e\le2(2d-3)+3+3d=7d-3.
\]
If $L=3$, the three pairwise intersecting edges form a star
in $K_{2,3}$. Each edge outside that star is disjoint from
exactly two of its edges. Hence at least two distinct witnesses
are needed, and
\[
\sum_{e\in\mathcal E}p_e\le3(2d-3)+2\cdot3+d=7d-3.
\]
In every case, \eqref{eq:max-k4-five-point-count} gives
\[
|\mathcal F|\le4d+5+7d-3=11n-53<12n-61,
\]
contrary to the assumed lower bound.

Thus $m\in\{2,3\}$. By~\eqref{eq:max-k4-full-cross},
$\mathcal F(23)$ consists of the $2n-9$ pairs meeting $\{4,5\}$,
together with one additional pair if $m=2$. Hence
$|\mathcal F(1)|+|\mathcal F(23)|=2n-6$, as required.
\end{proof}

\begin{lem}\label{lem:completion-pair-bounds}\label{clm01}
For $k\ge5$,
\[
\bigl(|\mathcal F(i)|+|\mathcal F([3]\setminus\{i\})|\bigr)\le
\begin{cases}
K,&|X_i|=k-1,\\
K-\delta,&2\le|X_i|\le k-2,\\
K+\delta,&|X_i|=1.
\end{cases}
\]
For $k=4$, write $N=n-3\ge6$. Then $K=3N-5$, $\delta=N-5$,
and the first and third bounds remain valid.
If $|\mathcal F|\ge g_2$ and $|X_i|=2$, then $\bigl(|\mathcal F(i)|+|\mathcal F([3]\setminus\{i\})|\bigr)=2N=K-\delta$.
For every $k\ge4$, if $X_i$ is a singleton and
\begin{equation}\label{eq:completion-small-core}
|\mathcal F(i)|<L:=\binom{n-4}{k-2}-\binom{n-k-2}{k-2},
\end{equation}
then $\bigl(|\mathcal F(i)|+|\mathcal F([3]\setminus\{i\})|\bigr)<K$.
\end{lem}
\begin{proof}
If $|X_i|=k-1$, then $\mathcal F(i)$ is a singleton and $\bigl(|\mathcal F(i)|+|\mathcal F([3]\setminus\{i\})|\bigr)=K$.
Let $2\le a=|X_i|\le k-2$ and $k\ge5$.
If a member of $\mathcal F([3]\setminus\{i\})$ avoids $X_i$, apply
Theorem~\ref{zh1} to $\mathcal F(i)$ and $\mathcal F([3]\setminus\{i\})$.
Its hypotheses hold since $n-3>(k-1)+(k-2)$, $k-2\ge3$,
and $a>1$. We obtain
\[
\bigl(|\mathcal F(i)|+|\mathcal F([3]\setminus\{i\})|\bigr)\le\binom{n-3}{k-2}-2\binom{n-k-2}{k-2}
+\binom{n-2k-1+a}{k-2}+2\le K-\delta.
\]
Otherwise
$\bigl(|\mathcal F(i)|+|\mathcal F([3]\setminus\{i\})|\bigr)\le Q_a:=\binom{n-3-a}{k-1-a}
+\binom{n-3}{k-2}-\binom{n-3-a}{k-2}$.
For $2\le a\le k-3$, Pascal's identity gives
$Q_{a+1}-Q_a=\binom{n-4-a}{k-3}-\binom{n-4-a}{k-1-a}\ge0$.
Indeed, $k-1-a\le k-3$ and
$(k-1-a)+(k-3)<n-4-a$, so unimodality applies.
Thus $Q_a$ is nondecreasing for $2\le a\le k-2$. Writing $z=n-k-3\ge k-2$ gives
$Q_{k-2}-(K-\delta)=z-\binom z{k-4}\le0$.
The assertion for $k=4$, $|X_i|=2$, is precisely
Lemma~\ref{lem:maximum-k4-core-two}.

If $X_i=\{x\}$, all $(k-2)$-sets through $x$ belong to $\mathcal F([3]\setminus\{i\})$.
When $\mathcal F([3]\setminus\{i\})(\overline{x})$ is nonempty, Theorem~\ref{ft}, applied
after deleting $x$ from $\mathcal F(i)$, gives $\bigl(|\mathcal F(i)|+|\mathcal F([3]\setminus\{i\})|\bigr)\le K$.
Otherwise a member of another single section avoids $x$, because
$\{i,x\}$ is not a cover. Counting the sets meeting this member
gives $\bigl(|\mathcal F(i)|+|\mathcal F([3]\setminus\{i\})|\bigr)\le K+\delta$.

To prove the final assertion, the residual family obtained by deleting
$x$ from $\mathcal F(i)$ is nontrivial. If $k\ge5$, the strict
Frankl--Tokushige inequality applies when the opposite family has
at least two members; the empty and singleton cases follow directly
from~\eqref{eq:completion-small-core}.
For $k=4$, put $M=n-4\ge5$ and view these residual families as graphs
$G,J$. We need $|G|+|J|<2M-2$ when $|G|<2M-3$ and $G$ is nontrivial.
For $|J|\le1$ this is immediate; for nontrivial $J$ it follows from
Lemma~\ref{lem:completion-graphs}(i). If $J$ is a star, it has at most
two edges, since three star edges would make $G$ trivial. With two
edges, $|G|+|J|\le M+2<2M-2$. This proves the final assertion.
\end{proof}

\subsection{Local replacements and covering number}

For $E\in\mathcal F(i)$, let $\mathcal D(\mathcal F,i,E)$ be
obtained by replacing $\mathcal F(i)$ with $\{E\}$ and
$\mathcal F([3]\setminus\{i\})$ with
\[
\left\{B\in\binom{[4,n]}{k-2}:B\cap E\ne\emptyset\right\},
\]
while leaving all other exact sections unchanged. The lifted
sections are understood in this definition. In particular, the
new complementary section has uniformity $k-2$.

\begin{lem}\label{lem:completion-cover-test}\label{clm3}
Let $\{i,j,\ell\}=[3]$ and $z\in[4,n]$.
Then $\{i,j,z\}$ covers $\mathcal F$ if and only if $z\in X_\ell$.
Moreover, if $z,w\in[4,n]$ and $\{i,z,w\}$ covers $\mathcal F$,
then $X_i\subseteq\{z,w\}$.
Consequently, if $\mathcal F$ is critically shifted,
$z\notin X_j\cup X_\ell$, and either $|X_i|\ge3$ or
$|X_i|=2$ with $z\notin X_i$, then $s_{iz}(\mathcal F)=\mathcal F$.
\end{lem}
\begin{proof}
The first assertion follows from the definition of $X_\ell$.
For the second, maximality puts every $(k-2)$-set meeting $X_i$
in $\mathcal F([3]\setminus\{i\})$, so a point of $X_i\setminus\{z,w\}$ would
produce a member avoiding $\{i,z,w\}$.
Under the final hypotheses, no three-cover contains $\{i,z\}$;
Proposition~\ref{p2} gives the conclusion.
\end{proof}

\begin{lem}\label{lem:completion-hub}
Every three-cover of $\mathcal G_1(n,k)$ contains its hub.
Relative to such a cover, every single section at a nonhub point
has a core of size at least $k-2$.
\end{lem}
\begin{proof}
Write the exceptional members as $S+a,S+b,T$, where
$|S|=k-1$, $|T|=k$, $S\cap T=\emptyset$, and $a,b\in T$ are
distinct. Their two-covers are
the edges of $K_{S,T}\cup\{ab\}$, whose vertex-cover number is $k$.
Here $K_{S,T}$ denotes the complete bipartite graph with parts
$S,T$. A three-cover avoiding the hub would cover every edge of
this graph: otherwise a disjoint edge extends, together with the
hub, to a member of $\mathcal G_1$ avoiding the proposed cover.
Thus every three-cover contains the hub and is completed by either
$\{a,b\}$ or $\{s,t\}$ with $s\in S$, $t\in T$.
The corresponding nonhub single sections are singletons, except
when $t\notin\{a,b\}$, in which case their cores have sizes
$k-2$ and $k-1$.
\end{proof}

\begin{lem}\label{lem:completion-collapse}
For $E\in\mathcal F(i)$, replace $\mathcal F(i)$ by $\{E\}$ and
$\mathcal F([3]\setminus\{i\})$ by all $(k-2)$-sets on $[4,n]$ meeting $E$, leaving
the other sections unchanged. Put
$\mathcal D=\mathcal D(\mathcal F,i,E)$. Then $\mathcal D$ is intersecting,
$|\mathcal D|-|\mathcal F|=K-\bigl(|\mathcal F(i)|+|\mathcal F([3]\setminus\{i\})|\bigr)$, and
\[
\tau(\mathcal D)=3
\quad\Longleftrightarrow\quad
E\cap X_j\cap Y_\ell=\emptyset
\quad\text{for both choices }\{i,j,\ell\}=[3].
\]
Every intersecting completion of $\mathcal D$ retains $[3]$ as a cover.
\end{lem}
\begin{proof}
Intersection and the size formula follow directly from the section
description. Since the full section $\mathcal F(123)$ and all three nonempty
single sections are retained, every two-cover has the form
$\{\ell,z\}$ with $\ell\in[3]$ and $z\in[4,n]$.
The enlarged $\mathcal F([3]\setminus\{i\})$ has empty core, excluding $\ell=i$;
for $\ell\ne i$, such a pair is a cover exactly when
$z\in E\cap X_j\cap Y_\ell$.
Finally, the full section $\mathcal F(123)$ prevents any completion from adding
a member avoiding $[3]$.
\end{proof}

\subsection{Intermediate cores}

\begin{lem}\label{lem:completion-intermediate-collapse}
Assume $\mathcal F$ is critically shifted, $|\mathcal F|\ge g_2$,
$X_1=\emptyset$, and $1\le|X_2|\le|X_3|$.
If either $X_2$ or $X_3$ has size between $2$ and $k-2$, a local replacement
at such a section preserves covering number three and increases
the size by at least $\delta$.
\end{lem}
\begin{proof}
The increase follows from Lemma~\ref{lem:completion-pair-bounds}.
It remains to choose a section and a retained member satisfying the
cover criterion in Lemma~\ref{lem:completion-collapse}.
The absence of a two-cover gives $X_2\cap X_3\cap Y_1=\emptyset$.
If $X_2=\{x\}$, apply a local replacement at $3$, choosing a member avoiding $x$
when $x\in Y_1$; then $x\notin X_3$, so this is possible.

Suppose $2\le|X_2|\le k-2$ and $|X_3|\ge3$.
For $z\in X_3\cap Y_1$, the absence of a two-cover gives
$z\notin X_2$. Since $X_1=\emptyset$ and $|X_3|\ge3$,
Lemma~\ref{lem:completion-cover-test} gives
$s_{3z}(\mathcal F)=\mathcal F$.
If a member of $\mathcal F(2)$ contained $z$, its shift would give
a member of $\mathcal F(23)$ avoiding $z$, a contradiction.
Every member of $\mathcal F(2)$ therefore satisfies the local replacement criterion.

The remaining case is $|X_2|=|X_3|=2$.
Put $Z=X_3\cap Y_1$, which is disjoint from $X_2$.
If $|Z|\le1$, choose a member of $\mathcal F(2)$ avoiding $Z$.
Otherwise $X_3\subseteq Y_1$ and $X_2\cap X_3=\emptyset$.
For $k\ge5$, suppose every member of $\mathcal F(2)$ meets $X_3$.
Choose $x\in X_3$ and $E\in\mathcal F(2)$ avoiding $x$;
such an $E$ exists because $x\notin X_2$.
Since $|E|=k-1\ge4$ and $|E\cap(X_2\cup X_3)|\le3$,
choose $z\in E\setminus(X_2\cup X_3)$.
Lemma~\ref{lem:completion-cover-test} gives
$s_{3z}(\mathcal F)=\mathcal F$.
Shifting $\{2\}\cup E$ shows that
$E\setminus\{z\}\in\mathcal F(23)$, contradicting
$x\in X_3\subseteq Y_1$.
Thus a member avoiding $X_3$ exists and may be retained in the local replacement.
For $k=4$, the failure to find a member avoiding $X_3$ instead
forces $\mathcal F(2)=\{X_2+z:z\in X_3\}$.
If $\mathcal F(23)\ne\emptyset$, then
$\mathcal F(23)=\{X_3\}$, and local replacement at $3$ is allowed because
$X_2\cap Y_1=\emptyset$.
If $\mathcal F(23)=\emptyset$, Lemma~\ref{lem:completion-cover-test}
makes $s_{2z}$ invariant for every $z\notin X_2\cup X_3$.
These shifts force
$\mathcal F(3)=\{X_3+z:z\in X_2\}$:
otherwise shifting such a member would produce a member of
$\mathcal F(23)$. Every member of $\mathcal F(1)$ must now contain
at least two points of $X_2\cup X_3$, so
$|\mathcal F(1)|\le6(N-4)+4=6N-20$, where $N=n-3$.
Both other pair sums equal $2N$, giving
$|\mathcal F|\le11N-20<12N-25=g_2$, a contradiction.
\end{proof}

\subsection{Singleton cores}

\begin{lem}\label{lem:completion-singleton-cases}
Assume $\mathcal F$ is critically shifted and $X_1=\emptyset$.
If $(|X_2|,|X_3|)=(1,1)$ or $(1,k-1)$, then either
$|\mathcal F|<g_2$, or there is a larger intersecting family
$\mathcal D$ with covering number three, containing the full
$123$-section and contained in no copy of $\mathcal G_1$.
\end{lem}
\begin{proof}
		We consider the two possible core configurations separately.
		
		\smallskip
		\noindent\textbf{Case 1:} $|X_2|=|X_3|=1$.
		Write $X_2=\{x\}$ and $X_3=\{y\}$.
If $|\mathcal F(2)|<L$, apply a local replacement at $2$, retaining a member avoiding $y$ when
$y\in Y_1$ and $y\ne x$. When $x=y$, the absence of a
two-cover gives $x\notin Y_1$.
Lemma~\ref{lem:completion-pair-bounds} makes the increase positive.
The new family is contained in no copy of $\mathcal G_1$. Indeed,
the unchanged empty core at $1$ would force $1$ to be the hub,
whereas the unchanged singleton core at $3$ would contradict
Lemma~\ref{lem:completion-hub}. The same argument applies at $3$.
We may therefore assume $|\mathcal F(2)|,|\mathcal F(3)|\ge L$.

Suppose $k\ge5$ and $x\ne y$.
In each section, restrict to members avoiding the other core point
and then delete its own core point.
The two resulting cross-intersecting $(k-2)$-uniform families on
$n-5$ points both have size at least
$\binom{n-5}{k-2}-\binom{n-k-2}{k-2}$.
Kruskal--Katona bounds either by $n-k-2$, whereas the displayed
lower bound exceeds $n-k-2$ for $k\ge5$, a contradiction.
For the lexicographic bound, write $v=n-5$ and $s=k-2$. The initial
lexicographic segment
of size $\binom vs-\binom{v-s+1}s$ consists of all $s$-sets
meeting a fixed $(s-1)$-set. Its $s$-uniform cross-family consists
of the $v-s+1=n-k-2$ sets containing that fixed set.
Thus $x=y$.
The same lexicographic computation, using $\mathcal F(2,x)$,
gives $0<|\mathcal F(1,\overline{x})|\le n-k-2$.
Here the initial segment of size $L$ consists of all $(k-2)$-sets
meeting a fixed $(k-2)$-set on $n-4$ points, and its
$(k-1)$-uniform cross-family has size $n-k-2$.
Also $\mathcal F(23,\overline{x})\ne\emptyset$, since $\{1,x\}$
is not a cover. Put $s=k-2$ and $v=n-4$.
Lemma~\ref{lem:completion-small-family}, applied to
$\mathcal F(1,\overline{x})$ and $\mathcal F(23,x)$, gives
\[
|\mathcal F(1,\overline{x})|+|\mathcal F(23,x)|
\le\binom v{s-1}-\binom{n-k-3}{s-1}+1.
\]
If the second family is empty, the same bound follows from
$|\mathcal F(1,\overline{x})|\le n-k-2$.
If $\mathcal F(1,x)\ne\emptyset$, Theorem~\ref{ft} gives
$|\mathcal F(1,x)|+|\mathcal F(23,\overline{x})|
\le\binom vs-\binom{v-s}s+1$.
Since $\binom{v-s-1}{s-1}\ge1$ and $v-s-1=n-k-3$, this implies
\[
|\mathcal F(1,x)|+|\mathcal F(23,\overline{x})|
\le\binom vs-\binom{n-k-3}s.
\]
When $\mathcal F(1,x)$ is empty, use any member of
$\mathcal F(1,\overline{x})$ to obtain the latter inequality.
Addition gives $\bigl(|\mathcal F(1)|+|\mathcal F(23)|\bigr)\le K$, and Lemma~\ref{lem:completion-strong-number}
gives $|\mathcal F|<g_2$.

Now suppose $k=4$ and write $N=n-3\ge6$.
If $x\ne y$, put
$G=\{A\setminus\{x\}:A\in\mathcal F(2),\ y\notin A\}$ and
$H=\{A\setminus\{y\}:A\in\mathcal F(3),\ x\notin A\}$.
These are cross-intersecting graphs on $N-2$ vertices, each with
at least $(2N-5)-(N-2)=N-3$ edges.
If $N\ge7$ and both are nontrivial, Lemma~\ref{lem:completion-graphs}(i)
would give $|G|+|H|\le6<2(N-3)$.
Thus one is a star. It has at least three edges, so every edge
of the other graph contains its center. The lower bounds force
both graphs to be the same full star.
For $N=6$, the four-point case of that lemma gives
$|G|+|H|\le6$, so $|G|=|H|=3$.
If either is a star, the same argument gives a common full star;
otherwise both have empty core.

In this latter case, each original section has size exactly seven:
it has three members avoiding the other core point, at most four
containing both core points, and at least seven in total.
The graph $\mathcal F(2,x)$ therefore contains the full four-edge
star at $y$ together with the nontrivial three-edge graph $G$.
It has no two-vertex cover: such a cover would have to contain
$y$, and its other vertex would then belong to every edge of $G$.
Consequently $\mathcal F(13,\overline{x})=\emptyset$.
Maximality gives the full $x$-star as $\mathcal F(13)$; similarly,
$\mathcal F(12)$ is the full $y$-star. Hence $\bigl(|\mathcal F(2)|+|\mathcal F(13)|\bigr)=\bigl(|\mathcal F(3)|+|\mathcal F(12)|\bigr)=7+5=12=2N$.
A local replacement at $2$, with the retained member avoiding $y$ when
necessary, increases the size by $\delta$. The unchanged cores
exclude containment in a copy of $\mathcal G_1$, as above.
In the common-star case, let $z$ be its center.
Each section has exactly $N-3$ members avoiding the other core
point and must contain all $N-2$ triples through both $x,y$.
Thus the two sections are precisely all triples through $x$
meeting $\{y,z\}$ and all triples through $y$ meeting $\{x,z\}$.
Their opposite graphs consist of the full star at the corresponding
core point together with the pair of the other two points.
Their pair sums are therefore $(2N-5)+N=K$.
Cross-intersection now forces every member of $\mathcal F(1)$ to
contain at least two points of $\{x,y,z\}$.
It follows that $|\mathcal F(1)|\le3(N-3)+1=3N-8$.
Lemma~\ref{lem:completion-graphs}(ii) gives $\bigl(|\mathcal F(1)|+|\mathcal F(23)|\bigr)\le K$, and
$|\mathcal F|\le N+3K<g_2$.

If $x=y$, the lexicographic bound gives
$0<|\mathcal V|\le N-3$, where $\mathcal V=\mathcal F(1,\overline{x})$.
Write $J=\mathcal F(1,x)$ and $P=\mathcal F(23,\overline{x})\ne\emptyset$.
If $J=\emptyset$, Lemma~\ref{lem:completion-graphs}(ii) gives
$\bigl(|\mathcal F(1)|+|\mathcal F(23)|\bigr)<K$. Otherwise Lemma~\ref{lem:completion-graphs}(i), on $N-1$ points,
gives $|J|+|P|\le2N-4$.
Maximality identifies $\mathcal F(23,x)$ with
$\{\{z\}:z\in\bigcap\mathcal V\}$.
Since $1\le|\mathcal V|\le N-3$, we have
$|\mathcal V|+|\bigcap\mathcal V|\le N-1$.
Consequently $\bigl(|\mathcal F(1)|+|\mathcal F(23)|\bigr)\le(2N-4)+(N-1)=K$, and hence
$|\mathcal F|\le N+K+2(K+\delta)=g_2$.
Equality would require
$\mathcal F(13)=\mathcal F(12)=\{\{x,z\}:z\ne x\}$.
For $z\notin\{x\}\cup\bigcap\mathcal V$, no three-cover contains
$\{2,z\}$. A third point $1$ or $3$ is excluded by
$X_3=\{x\}$ and $X_1=\emptyset$, respectively.
A third point outside $[3]$ would have to be $x$ to meet every
member of the full $x$-star $\mathcal F(13)$; but then a member
of $\mathcal V$ avoiding $z$ would be missed.
Hence $s_{2z}(\mathcal F)=\mathcal F$.
If $E\in\mathcal V$ contained such a $z$, shifting
$\{1\}\cup E$ would give
$E\setminus\{z\}\in\mathcal F(12)$, impossible because this pair
avoids $x$. Therefore every $E\in\mathcal V$ is contained in
$\bigcap\mathcal V$, and $\mathcal V$ is a singleton. Now $\bigl(|\mathcal F(1)|+|\mathcal F(23)|\bigr)\le(2N-4)+4=2N<K$, a contradiction.

\medskip\noindent\textbf{Case 2:} $|X_2|=1$ and $|X_3|=k-1$.
Write $X_2=\{x\}$ and $\mathcal F(3)=\{S\}$.
Since $|S|=k-1\ge3$, Lemma~\ref{lem:completion-cover-test}
gives $s_{3z}(\mathcal F)=\mathcal F$ for every $z\ne x$.
Also $x\notin S\cap Y_1$, since otherwise $\{1,x\}$ would cover
$\mathcal F$. If $E\in\mathcal F(2)$ contained a point
$z\in S\cap Y_1$, invariance would give
$E\setminus\{z\}\in\mathcal F(23)$, contrary to $z\in Y_1$.
Thus every $E\in\mathcal F(2)$ avoids $S\cap Y_1$ and satisfies
the local replacement criterion at $2$.
For any such $E$, shifting at each $z\in E\setminus\{x\}$ gives
$k-2$ distinct members $E\setminus\{z\}$ of $\mathcal F(23)$,
all containing $x$. In particular, $|\mathcal F(23)|\ge k-2$
and $\mathcal F(23,x)\ne\emptyset$.
If $|\mathcal F(2)|<L$, the local replacement is positive and cannot have a
$\mathcal G_1$ completion: hub $1$ is forced, but there are at least
$(k-2)+1+1\ge4$ nonhub members.
For $k\ge5$ and $|\mathcal F(2)|\ge L$, the same lexicographic
calculation gives $0<|\mathcal F(1,\overline{x})|\le n-k-2$.
Apply Lemma~\ref{lem:completion-small-family} to
$\mathcal F(1,\overline{x}),\mathcal F(23,x)$ for the first cross-sum
estimate above. If $\mathcal F(23,\overline{x})\ne\emptyset$, the
second estimate follows exactly as before, including the case
$\mathcal F(1,x)=\emptyset$.
If $\mathcal F(23,\overline{x})=\emptyset$, the absence of the two-cover
$\{1,x\}$ forces $x\notin S$. Every member of $\mathcal F(1,x)$
then meets the $(k-1)$-set $S$, so its single neighborhood gives
$|\mathcal F(1,x)|\le\binom{n-4}{k-2}-\binom{n-k-3}{k-2}$,
which is again the second estimate.
Adding gives $\bigl(|\mathcal F(1)|+|\mathcal F(23)|\bigr)\le K$ and then $|\mathcal F|<g_2$
by Lemma~\ref{lem:completion-strong-number}.
For $k=4$, choose $E\in\mathcal F(1)$ avoiding $x$.
Its three pairs belong to $\mathcal F(13,\overline{x})$ by the invariant
shifts. The nontrivial graph $\mathcal F(2,x)$ must therefore be
the triangle on $E$. Thus $|\mathcal F(2)|=3$, contradicting the standing assumption
$|\mathcal F(2)|\ge L=2N-5\ge7$.
Thus this last case cannot occur.
\end{proof}

\begin{lem}\label{lem:completion-all-cores}
If $\mathcal F$ is maximal intersecting, $\tau(\mathcal F)=3$,
$k=4$, and all three $X_i$ are nonempty, then
$|\mathcal F|\le g_2(n,4)$, with equality exactly when
$\mathcal F\cong\mathcal F(n,4)$.
\end{lem}
\begin{proof}
It suffices to consider $|\mathcal F|\ge g_2$. Write $N=n-3\ge6$.
By Lemma~\ref{lem:completion-pair-bounds}, the pair sum is at most $Q=4N-10$ for a singleton core, equals
$K=3N-5$ for a full core, and equals $2N<K$ for a two-point core.
If there are at most two singleton cores, then
$|\mathcal F|\le N+K+2Q=g_2$.

If all three cores are singletons and some single section has size
less than $2N-5$, its pair sum is less than $K$, giving strictness.
Otherwise all three single sections have size at least $2N-5$.
If two core points are different, the residual graph argument in
Lemma~\ref{lem:completion-singleton-cases} bounds both corresponding
pair sums by $K$, including the $N=6$ case. Hence
$|\mathcal F|\le N+2K+Q<g_2$.
If all three core points equal $x$, every $\mathcal F([3]\setminus\{i\})$ has a
member avoiding $x$, since $\{i,x\}$ is not a cover.
The graph cross-sum bound gives $\bigl(|\mathcal F(i)|+|\mathcal F([3]\setminus\{i\})|\bigr)\le K$ for all $i$, again
giving strictness.

Equality therefore requires precisely two singleton cores, say
$X_1=\{x\}$ and $X_2=\{y\}$, with $\bigl(|\mathcal F(1)|+|\mathcal F(23)|\bigr)=\bigl(|\mathcal F(2)|+|\mathcal F(13)|\bigr)=Q$, and a full
third core $\mathcal F(3)=\{E\}$ with $\bigl(|\mathcal F(3)|+|\mathcal F(12)|\bigr)=K$.
A singleton core can have pair sum $Q>K$ only if the complementary
graph is the full star at the core point.
Choosing a member $E_0$ of another single section avoiding $x$,
equality in the neighborhood count forces
\[
\mathcal F(1)=\{\{x\}\cup e:e\in\binom{[4,n]\setminus\{x\}}2,
\ e\cap E_0\ne\emptyset\}.
\]
The only three-set avoiding $x$ and meeting every member of this
family is $E_0$: the graph of all pairs meeting $E_0$ has $E_0$
as its unique three-vertex cover when its ground set has at least
five points. Thus another single section has at most one member
avoiding $x$.
If $y\ne x$, however, $\mathcal F(2)$ has $3N-9$ members and at
most $N-2$ of them contain $x$, a contradiction. Hence $x=y$.
The absence of the two-cover $\{1,x\}$ now gives $x\notin E$.
Consequently both $\mathcal F(1),\mathcal F(2)$ are exactly the
triples through $x$ meeting $E$, their opposite graphs are full
$x$-stars, and $\mathcal F(12)$ consists of all pairs meeting $E$.
Together with the full section $\mathcal F(123)$, this is the construction
$\mathcal F(n,4)$ defined in the introduction.
\end{proof}

\subsection{Empty cores}

\begin{lem}\label{lem:completion-two-empty}
Assume $\mathcal F$ is critically shifted,
$X_1=X_2=\emptyset$, $|X_3|\in\{0,1,k-1\}$, and $\bigl(|\mathcal F(1)|+|\mathcal F(23)|\bigr)<K$.
There is an intersecting family $\mathcal D$ with
$\tau(\mathcal D)=3$, $|\mathcal D|>|\mathcal F|$, containing the
full $123$-section and contained in no copy of $\mathcal G_1$.
\end{lem}
\begin{proof}
Apply the local replacement at $1$, choosing $E$ avoiding $X_3$ when $|X_3|\le1$.
The criterion is immediate.
If $|X_3|=k-1$, then $\mathcal F(3)=\{X_3\}$ and
$\mathcal F(12)$ consists of all $(k-2)$-sets meeting $X_3$.
By Lemma~\ref{lem:completion-cover-test}, $s_{3z}$ is invariant
for every $z\in[4,n]$.
For $E\in\mathcal F(1)$ and $z\in E$, shifting
$\{1\}\cup E$ therefore gives
$E\setminus\{z\}\in\mathcal F(13)$.
Hence $\mathcal F(13)$ contains the entire $(k-2)$-shadow of each
$E\in\mathcal F(1)$. In particular, $Y_2=\emptyset$ and
$|\mathcal F(13)|\ge k-1$, so any local replacement at $1$ satisfies the
criterion.

If $\mathcal D$ were contained in $\mathcal G\cong\mathcal G_1$,
then $[3]$ would cover $\mathcal G$.
The retained empty core $X_2$ and Lemma~\ref{lem:completion-hub}
force its hub to be $2$.
If $|X_3|\le1$, the retained section at $3$ contradicts that lemma.
Otherwise the single sections at $1,3$ and the retained
$|\mathcal F(13)|\ge k-1$ members give more than three members
avoiding this hub, also impossible.
\end{proof}

\subsection{Two full cores}

\begin{lem}\label{lem:completion-two-full-cores}
Suppose $X_1=\emptyset$, $\mathcal F(2)=\{C\}$,
$\mathcal F(3)=\{D\}$, and $|\mathcal F(23)|\ge2$.
Then $|\mathcal F|<g_2(n,k)$.
\end{lem}
\begin{proof}
All neighborhoods below are on $[4,n]$.
Put $r=k-1$, $a=n-k-2$, $b=n-2k$, and write
$\mathcal B=\mathcal F(23)$, $m=|\mathcal B|$,
$Y=\bigcap\mathcal B$, $y=|Y|\le k-3$.
Since $\{1,z\}$ is not a cover, $Y\cap C\cap D=\emptyset$.
Maximality gives
\begin{equation}\label{eq:completion-neighborhood}
|\mathcal F|=\binom{n-1}r-2\binom a{r-1}+2
-\bigl(|N_r(\mathcal B\cup\{C,D\})|-m\bigr).
\end{equation}

Suppose first that a member of $\mathcal F(1)$ avoids $Y$.
Let $\mathcal C=\binom{[4,n]}r\setminus N_r(\mathcal B)$ be the
full cross-family of $\mathcal B$. Since $\mathcal F(1)\subseteq\mathcal C$,
some member of $\mathcal C$ avoids $Y$.
Apply Theorem~\ref{zh1} to $\mathcal B$ and $\mathcal C$, with
uniformities $r-1,r$ and common core $Y$.
Its hypotheses hold because $r\ge3$, $n-3>2r-1$, and $y\ge0>-1$.
Using $|\mathcal C|=\binom{n-3}r-|N_r(\mathcal B)|$ gives
\[
|N_r(\mathcal B)|-m\ge2\binom{a+1}r-\binom{b+y+1}r-2.
\]
Write $s=|C\setminus D|=|D\setminus C|$,
$p=|Y\cap(C\setminus D)|$, $q=|Y\cap(D\setminus C)|$, and
$d=|Y\setminus(C\cup D)|$. Then $y=p+q+d$ and $s\ge p,q$.
Set $f(z)=\binom zr$. Because every member of $\mathcal B$ contains
$Y$, the new neighbors supplied by $C,D$ number at least
\[
G=2f(a)-f(a-s)-f(a-q-d)-f(a-p-d)+f(a-s-d).
\]
For $h\ge0$, write $\Delta_h(z)=f(z)-f(z-h)$.
Discrete convexity gives
\[
G-\Delta_y(a)=\Delta_p(a)-\Delta_p(a-q-d)
+\Delta_d(a-p)-\Delta_d(a-s)\ge0.
\]
Moreover $a\ge b+y+1$, so
$\Delta_y(a)\ge f(b+y+1)-f(b+1)$.
Consequently
\begin{equation}\label{eq:completion-neighborhood-lower}
|N_r(\mathcal B\cup\{C,D\})|-m
\ge2\binom{a+1}r-\binom{b+1}r-2.
\end{equation}

It remains that every member of $\mathcal F(1)$ meets $Y$.
Since $X_1=\emptyset$, this forces $y\ge2$ and therefore $k\ge5$.
Maximality then gives
$\mathcal B=\{B\in\binom{[4,n]}{r-1}:Y\subseteq B\}$.
Put $v=n-3$. The family $\mathcal F(1)$ consists of all $r$-sets
meeting $Y,C,D$, and
\begin{equation}\label{eq:completion-three-meet}
|\mathcal F(1)|\le\binom vr-\binom{v-y}r
-\binom{v-r}r+\binom{v-r-y}r.
\end{equation}
To see this directly, put $p=|Y\cap C|$, $q=|Y\cap D|$,
$h=|C\cap D|$, and assume $p\ge q$.
Then $p+q\le y$ and $h\le r-p$.
In inclusion--exclusion the dependence on $h$ is
$f(v-2r+h)-f(v-y-2r+p+q+h)$, an increasing function.
Replacing $h$ by $r-p$ cancels the term involving $q$.
Convexity bounds the remaining sum
$f(v-y-r+p)+f(v-r-p)$ by $f(v-y-r)+f(v-r)$,
proving~\eqref{eq:completion-three-meet}.

The left side of~\eqref{eq:completion-neighborhood-lower} is therefore
at least
\[
L(y)=\binom{v-y}r-\binom{v-y}{r-1-y}
+\binom ar-\binom{a-y}r.
\]
We use the elementary inequality
$\binom MR-\binom{M-B}R\ge\binom Mq$ for
$B>R>q\ge0$ and $M\ge B+q$. For completeness, when
$M=B+q$ this is binomial unimodality, and when $q=0$ it is
immediate. For $q\ge1$, increasing $M$ by one and applying
Pascal's identity reduces the assertion to the same inequality
with $(R,q)$ replaced by $(R-1,q-1)$; induction proves it.
Taking $B=r$, $R=r-1$, $M=v-y-1$, $q=r-1-y$ gives
$L(y)\ge L(y+1)$ for $2\le y<r-2$.
Finally,
\[
L(r-2)-\left(2\binom{a+1}r-\binom{b+1}r-2\right)
=\binom a{r-2}-a>0,
\]
because $a\ge r\ge4$.
This proves~\eqref{eq:completion-neighborhood-lower} in all cases.
Substituting into~\eqref{eq:completion-neighborhood} gives an upper
bound that is smaller than $g_2$ by
$\binom{n-k-3}{k-2}-\binom{n-2k}{k-2}>0$.
\end{proof}

\subsection{Final structural alternative}

\begin{prop}[Critical-family alternative]\label{prop:completion-alternative}\label{123}
Let $n>2k$, $k\ge4$, and let $\mathcal F$ be maximal intersecting,
critically shifted, and satisfy $|\mathcal F|\ge g_2$.
After choosing a critically shifted labeling with cover $[3]$,
at least one of the following holds:
\begin{enumerate}
\item[\rm(i)] Some point is avoided by exactly three members.
\item[\rm(ii)] $k=4$ and $\mathcal F\cong\mathcal F(n,4)$.
\item[\rm(iii)] A local replacement of a section with
$2\le|X_i|\le k-2$ preserves covering number three and increases
the size by at least $\delta$. At least one unchanged single
section has empty core.
\item[\rm(iv)] There is a larger intersecting family $\mathcal D$
with covering number three, containing the full section $\mathcal F(123)$ and
contained in no copy of $\mathcal G_1$.
\end{enumerate}
\end{prop}
\begin{proof}
The required labeling exists by
Lemma~\ref{lem:verified-cover-labeling}.

\noindent\textbf{Case 1:} All three cores are nonempty. Lemmas
\ref{lem:completion-pair-bounds} and
\ref{lem:completion-strong-number} exclude $k\ge5$, while
Lemma~\ref{lem:completion-all-cores} gives (ii) for $k=4$.

\medskip\noindent\textbf{Case 2:} At least two cores are empty.
By symmetry, suppose $X_1=X_2=\emptyset$.
M\"ors's theorem gives, after interchanging $1,2$,
$|\mathcal F(1)|\le H_0:=\binom{n-4}{k-2}-\binom{n-k-3}{k-2}+1$.
At $n=2k+1$, where its strict ground-set hypothesis does not hold,
the same conclusion follows instead from
$|\mathcal F(1)|+|\mathcal F(2)|\le\binom{2k-2}{k-1}=2H_0$:
cross-intersecting families cannot contain complementary sets.
This sum bound also supplies the strict implication used below
at that boundary.
Lemma~\ref{lem:completion-small-section}, or
Lemma~\ref{lem:completion-graphs}(ii) if $k=4$, gives $\bigl(|\mathcal F(1)|+|\mathcal F(23)|\bigr)\le K$.
If $\bigl(|\mathcal F(1)|+|\mathcal F(23)|\bigr)=K$, its equality condition gives $|\mathcal F(1)|=H_0$;
M\"ors's strict implication then gives $|\mathcal F(2)|\le H_0$,
and hence $\bigl(|\mathcal F(2)|+|\mathcal F(13)|\bigr)\le K$.
If both pair sums equal $K$, the third is at most $K+\delta$
when its core is nonempty; when its core is empty, apply M\"ors
to $\mathcal F(1),\mathcal F(3)$ and the same small-section lemma.
This gives $|\mathcal F|<g_2$, by
Lemma~\ref{lem:completion-strong-number} for $k\ge5$, or by
$N+2K+(K+\delta)=11N-20<12N-25$ for $k=4$.
We may thus arrange $\bigl(|\mathcal F(1)|+|\mathcal F(23)|\bigr)<K$.
If $2\le|X_3|\le k-2$, any local replacement at $3$ satisfies the cover
criterion and gives (iii). Otherwise
Lemma~\ref{lem:completion-two-empty} gives (iv).

\medskip\noindent\textbf{Case 3:} Exactly one core is empty.
Relabel it as $X_1$ and
assume $1\le|X_2|\le|X_3|$.
An intermediate core gives (iii) by
Lemma~\ref{lem:completion-intermediate-collapse}.
The cases $(1,1)$ and $(1,k-1)$ give (iv) or a contradiction by
Lemma~\ref{lem:completion-singleton-cases}.
Otherwise $|X_2|=|X_3|=k-1$.
Here $\mathcal F(23)\ne\emptyset$: otherwise the two members
avoiding $1$ would meet, and $1$ together with a point in their
intersection would form a two-cover.
Lemma~\ref{lem:completion-two-full-cores} excludes
$|\mathcal F(23)|\ge2$. Thus $|\mathcal F(23)|=1$, and exactly three
members avoid $1$, giving (i).
\end{proof}

\section{Proof of the Main Results}

\subsection{Three exceptional members}
The following estimate completes the three-member case that occurs
after the core reductions. No critical-shiftedness assumption is
needed for the numerical bound.

\begin{lem}\label{lem:completion-diversity-three}
Suppose that $\mathcal F$ is maximal intersecting, $\tau(\mathcal F)\ge3$,
and precisely three members avoid a point $x$. Then either
$\mathcal F\cong\mathcal G_1(n,k)$ or $|\mathcal F|\le q(n,k)$.
\end{lem}
\begin{proof}
Write the three exceptional members as $E_1,E_2,E_3$.
Their total intersection is empty. By maximality, the remaining
members are exactly the $k$-sets through $x$ meeting every $E_i$.
Order the three pairwise intersection sizes as $a\ge b\ge c\ge1$.
Then $a+b\le k$. Set $f(t)=\binom{n-2k-1+t}{k-1}$, interpreting a
binomial coefficient with negative upper argument as zero.
Inclusion--exclusion gives
\[
|\mathcal F|=\binom{n-1}{k-1}-3\binom{n-k-1}{k-1}+3
+f(a)+f(b)+f(c)-f(a+b+c-k).
\]
If $a=k-1$, then $b=c=1$, giving $\mathcal G_1$.
Otherwise $a\le k-2$. The first differences of $f$ are nonnegative
and nondecreasing. Increasing $c$ to $b$, and then $a$ to $k-b$,
does not decrease the last four terms. For $b\ge2$ their sum is
therefore at most
$f(k-b)+f(b)\le f(k-2)+f(2)$, by discrete convexity.
For $b=1$ it is at most
$f(k-2)+2f(1)-f(0)\le f(k-2)+f(2)$.
Substitution gives $|\mathcal F|\le q(n,k)$, as required.
\end{proof}

\begin{lem}\label{lem:verified-three-exceptional}
Let $\mathcal F\subseteq\binom{[n]}k$ be maximal intersecting,
with $\tau(\mathcal F)\ge3$, and admitting a critically shifted
labeling. If $|\mathcal F(\overline{x})|=3$ for some $x$, then
\[
\mathcal F\cong\mathcal R_r(n,k)
\quad\text{for some }1\le r\le\lfloor k/2\rfloor.
\]
\end{lem}

\begin{proof}
Fix a critically shifted labeling and write
$\mathcal F(\overline{x})=\{E_1,E_2,E_3\}$.
Maximality implies that the members containing $x$ are exactly
the $k$-sets containing $x$ and meeting each $E_i$.
The exceptional members are pairwise intersecting and have empty
total intersection; otherwise a pair containing $x$ would cover
$\mathcal F$.

Put $U_i=E_i\setminus\bigcup_{j\ne i}E_j$.
We show that at most one $U_i$ is nonempty.
Otherwise, after reindexing the exceptional members, take
$u\in U_1$, $v\in U_2$ with $u<v$, and put
$E_2'=(E_2\setminus\{v\})\cup\{u\}$.
Since $u\notin E_3$, we have $E_2'\ne E_3$.
If $E_2'=E_1$, then $E_1\setminus\{u\}\subseteq E_2$;
the empty triple intersection and $u\notin E_3$ would imply
$E_1\cap E_3=\emptyset$. Thus $E_2'\ne E_1$ as well.
Consequently, $s_{uv}$ changes $\mathcal F$, and its exceptional
members are $E_1,E_2',E_3$, still with empty total intersection.

Let $P$ be any pair avoiding $x$.
Choose $z_i\in E_i\setminus(P\cup\{u,v\})$.
These choices are possible: at most three points are excluded
from $E_1,E_2$, and at most two from $E_3$.
Since $n-4\ge k$, extend $\{x,z_1,z_2,z_3\}$ to a $k$-set
$H$ avoiding $P\cup\{u,v\}$.
Then $H\in\mathcal F$ by maximality, and $H$ is unchanged by
$s_{uv}$.
Thus $P$ does not cover the shifted family.
Nor can a pair containing $x$ cover it, because the exceptional
members have empty total intersection.
This contradicts critical shiftedness.

Reindex so that $U_1=U_2=\emptyset$, and set
$A=E_1\cap E_2,\quad B=E_1\cap E_3,\quad C=E_2\cap E_3,\quad D=U_3$.
These sets are pairwise disjoint, and
$E_1=A\cup B,\quad E_2=A\cup C,\quad E_3=B\cup C\cup D$.
Uniformity gives $|A|=k-r$, $|B|=|C|=r$, and $|D|=k-2r$.
Pairwise intersection gives $r\ge1$, while $|D|\ge0$ gives
$r\le\lfloor k/2\rfloor$.
Maximality now identifies $\mathcal F$ with $\mathcal R_r(n,k)$.
\end{proof}

\subsection{The maximum case}

\begin{thm}[Largest size]\label{prop:verified-numerical-maximum}
Let $n>2k$ and $k\ge4$. Every intersecting family
$\mathcal F\subseteq\binom{[n]}k$ with $\tau(\mathcal F)\ge3$
satisfies
\[
|\mathcal F|\le
\binom{n-1}{k-1}-3\binom{n-k-1}{k-1}
+\binom{n-k-2}{k-1}+\binom{n-2k}{k-1}+3
=g_1(n,k).
\]
Consequently $f(n,k,3)=g_1(n,k)$, also when restricted to
maximal intersecting families.
\end{thm}
\begin{proof}
The family $\mathcal G_1(n,k)$ gives
$f(n,k,3)\ge g_1(n,k)>g_2(n,k)$.
Choose an intersecting family $\mathcal F$ of maximum cardinality
among all families with covering number at least three, and put
$s=|\mathcal F|$. Such a family is maximal intersecting.
Apply admissible shifts until none changes it; every intermediate
family is again maximum and hence maximal. Since $s\ge g_1>g_2$, Theorem~\ref{fws} and
\eqref{gk3} exclude a fully shifted terminal family.
Thus the terminal family $\mathcal H$ is critically shifted.
By Proposition~\ref{p2} and Lemma~
\ref{lem:verified-cover-labeling}, it has covering number three
and admits a critically shifted labeling with cover $[3]$.

Apply Proposition~\ref{prop:completion-alternative} to $\mathcal H$.
Alternatives (iii) and (iv) each produce an intersecting family
with covering number three and cardinality strictly greater than
$s$, contrary to the definition of $s$.
Alternative (ii) has cardinality $g_2<g_1\le s$, and is also
impossible. Therefore exactly three members of $\mathcal H$
avoid some point. Lemma~\ref{lem:completion-diversity-three}
now gives either $\mathcal H\cong\mathcal G_1$, or
$s=|\mathcal H|\le q<g_1$. The latter contradicts $s\ge g_1$.
Hence $s=g_1$.

The bound is attained by the maximal family $\mathcal G_1$,
proving both assertions.
\end{proof}

\subsection{Connectivity and inverse shifts}

To recover uniqueness and control the gap below the maximum, we must
also understand families whose shifted images lie inside
$\mathcal G_1$. The required size deficit follows from a connectivity
bound.

For a family $\mathcal V$ of sets of one fixed size, its
\emph{disjointness graph} has vertex set $\mathcal V$ and joins
two sets when they are disjoint. If two nonempty disjoint
subfamilies $\mathcal A,\mathcal B\subseteq\mathcal V$ are
cross-intersecting, no edge joins them. Hence
$\mathcal V\setminus(\mathcal A\cup\mathcal B)$ is a vertex
separator, and its size is at least the vertex connectivity of
the graph. In the inverse-shift argument, the size of this separator is exactly
the deficit from $g_1$. The next lemma supplies the required
connectivity bound for the graphs arising from $\mathcal G_1$.

\begin{lem}[Connectivity of a disjointness graph]\label{lem:completion-connectivity}
Let $r\ge3$, $m\ge2r+1$, and let $S,T$ be disjoint subsets of
an $m$-element ground set, with $|S|\ge2$, $|T|\ge3$.
Fix distinct $a,b\in T$ and put $L=K_{S,T}\cup\{ab\}$.
The disjointness graph on the $r$-sets containing an edge of $L$,
after its isolated vertices are removed, has vertex connectivity
at least $\binom{m-r-2}{r-2}$.
\end{lem}
\begin{proof}
For disjoint edges $e,f$ of $L$, take the block consisting of the
$r$-sets containing $e$ and avoiding $f$, together with those
containing $f$ and avoiding $e$.
Its induced graph is the bipartite disjointness graph on two copies
of $\binom{[m]\setminus(e\cup f)}{r-2}$.
It is regular of degree $d=\binom{m-r-2}{r-2}$ and edge-transitive.
It is connected: two $(r-2)$-sets differing in one point have a
common disjoint $(r-2)$-set because $m-4\ge2(r-2)+1$, and such
exchanges connect each part.
Watkins's theorem~\cite{Watkins} states that a connected finite
edge-transitive graph has vertex connectivity equal to its minimum
degree. Thus the vertex connectivity of this block is $d$.

Every nonisolated vertex belongs to a block: choose an edge of $L$
inside it and another edge inside a disjoint neighbor.
If $f,g$ are intersecting partner edges disjoint from $e$, the
blocks indexed by $\{e,f\}$ and $\{e,g\}$ overlap in at least
$\binom{m-5}{r-2}\ge d$ vertices.
We show that these overlaps join all the blocks into a connected
union.
Indeed, the disjointness graph on the edges of $L$ is connected:
on the cross edges it joins pairs in $S\times T$ whose two
coordinates differ, and the additional edge $ab$ has a disjoint
cross edge. For each fixed $e$, its disjoint partners form the
edge set of a connected graph, namely $L-V(e)$ with isolated
vertices omitted, so their line graph is connected.
Thus any two blocks can be joined by a sequence of blocks whose
consecutive overlaps have the stated size.
Deleting fewer than $d$ vertices leaves each block connected and
each consecutive overlap nonempty; their union remains connected.
\end{proof}

\begin{lem}[Inverse-shift gap]\label{lem:completion-inverse-shift}
Let $\mathcal F$ be maximal intersecting with $\tau(\mathcal F)\ge3$.
If $s_{ij}(\mathcal F)$ has covering number at least three and is
contained in a copy of $\mathcal G_1(n,k)$, then either
$\mathcal F\cong\mathcal G_1(n,k)$ or
\[
|\mathcal F|\le\max\{q(n,k),g_2(n,k)-1\}.
\]
\end{lem}
\begin{proof}
Write $\mathcal H=s_{ij}(\mathcal F)\subseteq\mathcal G$,
where $\mathcal G$ has hub $x$ and exceptional members
$B_1=S+a$, $B_2=S+b$, $B_3=T$.
All three exceptional members belong to $\mathcal H$; any proper
subfamily of them has nonempty intersection and would give a
two-cover through $x$.
If $x\notin\{i,j\}$, the shift preserves the number of members
avoiding $x$. If $x=j$, it can only increase that number.
Since at least three members avoid every point in a family of
covering number at least three, either case gives exactly three
members of $\mathcal F$ avoiding $x$.
Lemma~\ref{lem:completion-diversity-three} applies.

It remains that $x=i$; write $y=j$.
Neither $y=a$ nor $y=b$ is possible: an exceptional member through
$y$ would have its shifted mate outside $\mathcal G$, whereas
$s_{xy}(\mathcal H)=\mathcal H$.
On $U=[n]\setminus\{x,y\}$, let $\mathcal V$ consist of the
$(k-1)$-sets $R$ with $x+R\in\mathcal G$.
These are precisely the sets containing an edge of
$K_{S\setminus\{y\},T\setminus\{y\}}\cup\{ab\}$.
The two part sizes are $(k-2,k)$, $(k-1,k-1)$, or $(k-1,k)$,
according as $y\in S$, $y\in T\setminus\{a,b\}$, or $y\notin S\cup T$.
The isolated vertices of its disjointness graph are exactly
$\mathcal Z=\{B_\nu\setminus\{y\}:y\in B_\nu\}$.
Indeed, if an edge of the kernel graph avoids a $(k-1)$-set
$R$, it extends to a $(k-1)$-set disjoint from $R$, because
$|U\setminus R|=n-k-1\ge k$. Thus $R$ is isolated exactly when
it covers that graph.
A cover of a complete bipartite graph contains an entire part.
If $y\in S$, the only covers of size $k-1$ are
$(S\setminus\{y\})+a$ and $(S\setminus\{y\})+b$.
If $y\in T\setminus\{a,b\}$, the only such cover is
$T\setminus\{y\}$: choosing the $S$-part would leave $ab$ uncovered.
If $y\notin S\cup T$, neither part yields a cover of size $k-1$.
These are precisely the members of $\mathcal Z$.

For $R\in\mathcal Z$, both $x+R$ and $y+R$ belong to
$\mathcal H$, and hence to $\mathcal F$.
The exceptional members avoiding $y$ are unchanged by the shift.
Thus all three $B_\nu$ belong to $\mathcal F$.
Every member of $\mathcal G$ containing both $x,y$ meets
$\mathcal F$, so maximality puts all of them in $\mathcal F$.
The remaining members are described by disjoint subfamilies
$\mathcal A,\mathcal B\subseteq\mathcal V\setminus\mathcal Z$,
according as $x+R$ or $y+R$ belongs to $\mathcal F$.
They are cross-intersecting, and
\[
g_1-|\mathcal F|
=|\mathcal D|,\qquad
\mathcal D=(\mathcal V\setminus\mathcal Z)
\setminus(\mathcal A\cup\mathcal B).
\]
If either $\mathcal A$ or $\mathcal B$ is empty, then
$\mathcal F$ is contained in $\mathcal G$ or its image under
$(x\ y)$. Maximality gives equality.
Otherwise $\mathcal D$ separates two nonempty parts of the
nonisolated disjointness graph.
Lemma~\ref{lem:completion-connectivity} gives
$g_1-|\mathcal F|\ge\binom{n-k-3}{k-3}=g_1-g_2+1$,
as required.
\end{proof}

\subsection{Uniqueness and the second largest size}

The inverse-shift bound lets us transfer the structure of a terminal
family back through the shift sequence. We first apply it at the maximum.

\begin{thm}[Uniqueness of the maximum]\label{thm:completion-unique-maximum}
For $n>2k$ and $k\ge4$, every intersecting family with covering
number at least three and size $g_1(n,k)$ is isomorphic to
$\mathcal G_1(n,k)$.
\end{thm}
\begin{proof}
Theorem~\ref{prop:verified-numerical-maximum} implies that such a
family is maximal intersecting.
Perform admissible shifts until none changes it. All intermediate
families remain maximal because their size is the numerical maximum.
The terminal family is not fully shifted, by Theorem~\ref{fws}
and~\eqref{gk3}, and is therefore critically shifted.
Normalize its minimum cover using
Lemma~\ref{lem:verified-cover-labeling}.
In Proposition~\ref{prop:completion-alternative}, both alternatives
producing a larger family contradict the numerical maximum,
and (ii) has smaller size. Thus (i) holds, and
Lemma~\ref{lem:completion-diversity-three}, with $q<g_1$, identifies
the terminal family as $\mathcal G_1$.
Lemma~\ref{lem:completion-inverse-shift} propagates this conclusion
backward through the shift sequence, since each predecessor has
size $g_1>\max\{q,g_2-1\}$.
\end{proof}

\begin{proof}[Proof of Theorem~\ref{thm:second-size}]
The constructions in the introduction establish sharpness. It remains to prove the upper bound.
Suppose that a maximal intersecting family with covering number
at least three, not isomorphic to $\mathcal G_1$, has size greater than
$M=\max\{q,g_2\}$. Among all such families choose one of largest
size $s$. By Theorem~\ref{thm:completion-unique-maximum}, $s<g_1$.

Apply a changing shift that preserves covering number at least
three, and complete its image to a maximal intersecting family
$\mathcal J$. Adding members cannot decrease the covering number, so
$\tau(\mathcal J)\ge3$.
If $\mathcal J\cong\mathcal G_1$, the inverse-shift gap contradicts
$s>M$. Otherwise the choice of $s$ forces $|\mathcal J|=s$,
so the shifted image was already maximal and may replace
$\mathcal F$. Each changing shift strictly decreases
$\sum_{F\in\mathcal F}\sum_{x\in F}x$, so the procedure terminates. The terminal family cannot be fully shifted:
by~\eqref{gk3}, its cardinality satisfies
$s>M\ge g_2>|\mathcal K(n,k,3)|$, contrary to
Theorem~\ref{fws} if it were fully shifted.
It is therefore critically shifted and has the same size $s$.

Apply Proposition~\ref{prop:completion-alternative}.
Alternative (i) contradicts
Lemma~\ref{lem:completion-diversity-three}, while (ii) has size $g_2<s$. Alternative (iii) gives a family of size at least
$s+\delta>g_2+\delta=g_1$, contradicting the numerical maximum.
In alternative (iv), a maximal completion has covering number
at least three and is not $\mathcal G_1$, but has size larger than $s$,
contradicting its choice. This proves the upper bound throughout
the full range $n>2k$.
\end{proof}

\subsection{Equality under critical shiftedness}

It remains to determine when a critically shifted family attains the
second largest size. The next lemma resolves the case in which a
local replacement produces a copy of $\mathcal G_1$.

\begin{lem}[Equality after a local replacement]\label{lem:completion-collapse-equality}
Let $n>2k$ and $k\ge4$, and let $\mathcal F\subseteq\binom{[n]}k$
be maximal intersecting with minimum cover $[3]$.
Suppose $X_1=\emptyset$, $2\le|X_2|\le k-2$, and a local replacement
at $2$ produces a copy of $\mathcal G_1(n,k)$.
Then $\mathcal F\cong\mathcal G_p(n,k)$ for some $2\le p\le k-1$.
If its size increase is $\delta$, then $p=2$, or $k=4$ and $p=3$.
\end{lem}
\begin{proof}
The unchanged empty core at $1$ forces the hub of the resulting
copy of $\mathcal G_1$ to be $1$.
Use its exceptional members $S+a,S+b,T$ as before.
The other two points of the cover are either $a,b$, or a pair
$s,t$ with $s\in S,t\in T$.

We use the following elementary duality observation. If a family consists of all
$(k-1)$-sets meeting each of at most three fixed sets of size at
most $k-1$, then every $(k-1)$-set cross-intersecting with it
contains one of those fixed sets. Otherwise choose a point of
each fixed set outside the proposed transversal and extend the
at most three chosen points to a disjoint $(k-1)$-set.
This extension is possible on $[4,n]$, since
$n-k-2\ge k-1$ and $k-1\ge3$.

If the cover pair is $a,b$, the unchanged section $\mathcal F(1)$
consists of all $(k-1)$-sets meeting $S$ and $T\setminus\{a,b\}$.
The dual observation forces every member of the original
$\mathcal F(2)$ to be $S$ or to contain $T\setminus\{a,b\}$.
Since it contains the retained member $S$ and at least one other
member, its core has size at most one, a contradiction.

If the pair is $s,t$ with $t\notin\{a,b\}$, the two off-hub
single sections of the new family have sizes two and one.
The replaced section must be the singleton $T\setminus\{t\}$.
The three fixed sets defining $\mathcal F(1)$ are this set and
$(S\setminus\{s\})+a,(S\setminus\{s\})+b$.
Besides the retained set $T\setminus\{t\}$, the only
possibilities are $(S\setminus\{s\})+a$ and
$(S\setminus\{s\})+b$, because all three defining sets have
size $k-1$. Each of the latter two meets the retained set in
exactly one point. The original section has at least two members,
so its core again has size at most one, a contradiction.

It remains that the cover pair is $s,a$, up to exchanging $a,b$.
Put $C=S\setminus\{s\}$ and $E=T\setminus\{a\}$.
The off-hub single sections are $\{C+b\}$ and $\{E\}$, and
the pair section avoiding the hub is $\{C\}$.
The unchanged $\mathcal F(1)$ consists of all $(k-1)$-sets
meeting $C,E$. Therefore the original section being replaced consists
of sets containing $C$ or equal to $E$.
A core of size at least two excludes any mixture of these two
types, and excludes collapsing the section $\{E\}$.
The original section is consequently
$\{C+z:z\in Z\}$ for some $Z\subseteq E$, $|Z|=p\ge2$;
cross-intersection with $\{E\}$ gives $Z\subseteq E$.
The cover point of this section is $s$, so the exceptional members
of the original family are precisely
$\{S+z:z\in Z\}\cup\{S+a,T\}$.
There are $p+1$ members of the form $S+z$, with distinct additional
points in the disjoint $k$-set $T$.
These are exactly the exceptional members defining $\mathcal G_p$.
Maximality then supplies every member through the hub meeting
them all, and determines the entire family.

The equality cases for a size increase of $\delta$ now follow
directly from~\eqref{eq:Gp-difference}.
\end{proof}

\begin{proof}[Proof of Theorem~\ref{mre}]
By Proposition~\ref{p2}, $\tau(\mathcal F)=3$,
so Theorem~\ref{thm:second-size} gives the upper bound.
Suppose $|\mathcal F|=M$. By
Lemma~\ref{lem:verified-cover-labeling}, choose a critically shifted
labeling in which $[3]$ is a cover.

Apply Proposition~\ref{prop:completion-alternative}.
Alternative (iv) is impossible: a maximal completion of its
larger family would retain covering number three, be contained
in no copy of $\mathcal G_1$, and contradict Theorem~\ref{thm:second-size}.
Alternative (ii) gives the listed family $\mathcal F(n,4)$.
In alternative (i), Lemma~\ref{lem:verified-three-exceptional}
gives $\mathcal F\cong\mathcal R_r$.
Here $r\ne1$, since $g_1>M$.
By~\eqref{eq:R-comparison}, equality holds precisely
when $r=2$ or $n=2k+1$, and the comparison of $q,g_2$
gives exactly parts (i) and (iii) of the theorem.

In alternative (iii), the numerical maximum implies
$M+\delta\le g_1=g_2+\delta$.
Thus $M=g_2$, the local replacement has gain exactly $\delta$, and its
image has size $g_1$.
By Theorem~\ref{thm:completion-unique-maximum}, that image is
$\mathcal G_1$. Two unchanged empty-core sections are impossible
by Lemma~\ref{lem:completion-hub}; otherwise relabel the unchanged
empty-core section as $1$ and the replaced section as $2$.
Lemma~\ref{lem:completion-collapse-equality} now gives
$\mathcal G_2$, or $\mathcal G_3$ when $k=4$.
The parameter restriction $M=g_2$ is exactly the one in part (ii).

Conversely, every listed family is maximal intersecting with covering
number three, admits a critically shifted labeling, and attains $M$
in the stated parameter range by the size formulas in the introduction. This completes the classification.
\end{proof}

\end{document}